\documentclass[reqno,12pt]{amsart}
\usepackage[margin=0.8in]{geometry}
\usepackage{amsmath}
\usepackage{amssymb}
\usepackage{dsfont}
\usepackage{graphicx} 
\usepackage{microtype}
\usepackage{amsthm}

\usepackage[dvipsnames]{xcolor}
\usepackage{hyperref}
\hypersetup{
    colorlinks=true,
    linkcolor=cyan!80!black,
    citecolor=MidnightBlue,
    urlcolor=magenta,
}
\usepackage{nicematrix}

\usepackage{tikz}
\tikzset{arrow/.style = {-stealth},}

\theoremstyle{plain}
\newtheorem{theorem}{Theorem}
\newtheorem{corollary}{Corollary}
\newtheorem{proposition}{Proposition}
\newtheorem{lemma}{Lemma}

\theoremstyle{definition}

\newtheorem{definition}{Definition}

\newcommand{\bfx}{\mathbf{x}}
\newcommand{\bfy}{\mathbf{y}}
\newcommand{\bbE}{\mathbb{E}}
\newcommand{\bbP}{\mathbb{P}}
\newcommand{\bbR}{\mathbb{R}}
\newcommand{\bbN}{\mathbb{N}}

\newcommand{\cN}{\mathcal{N}}

\newcommand{\cF}{\mathcal{F}}
\newcommand{\cG}{\mathcal{G}}
\newcommand{\cS}{\mathcal{S}}

\title{The two-dimensional $\ell_{\infty}$ Directed Spanning Forest}

\author[D. Pal]{Dipranjan Pal}
\address{
    Statistics and Mathematics Unit \\
    Indian Statistical Institute \\
    203 B.T. Road, Kolkata 700108 \\
    West Bengal, India.
}
\email{dipranjan\_r@isical.ac.in}

\author[K. Saha]{Kumarjit Saha}
\address{
    Department of Mathematics \\
    Ashoka University\\
    Plot no 2, Rajiv Gandhi Education City, Sonipat 131029\\
    Haryana, India.
}
\email{kumarjit.saha@ashoka.edu.in}

\begin{document}

\begin{abstract}
   We investigate the $\ell_{\infty}$ \textit{ directed spanning forest} (DSF), a directed forest whose vertex set is given by
  a homogeneous Poisson point process, in which each Poisson
  point connects to the nearest Poisson point, measured in $\ell_{\infty}$ distance, with a strictly larger $y$-
  coordinate. In this paper, we prove that the two-dimensional $\ell_{\infty}$ DSF is connected and obtain an optimal
  estimate on the tail decay of coalescing time of two $\ell_{\infty}$ DSF paths. Similar estimates were earlier obtained in \cite{MR4203342} for the $\ell_2$ (Euclidean) DSF and in \cite{garcia2025directed} for the $\ell_p$ DSF with $p \in [1, \infty]$ in two-dimension. The geometry of $\ell_{\infty}$ balls compels us to develop new arguments.
\end{abstract}

\maketitle

\section{Introduction and main results}\label{section_1}

Consider a homogeneous Poisson Point Process (PPP)
$\mathcal{N}$ with intensity $\lambda>0$ on $\mathbb{R}^2$ . For $\mathbf{x} \in \mathbb{R}^2$ let 
$\mathbf{x}(i)$ denote  the $i$-th coordinate of $\mathbf{x}$ for $i=1,2$. The $\ell_\infty$ directed
spanning forest (DSF) ancestor of $\mathbf{x}$ is denoted as $h(\mathbf{x}) \in \mathcal{N}$ and defined as
\begin{align}
\label{def:DSF_step}  
h(\mathbf{x}) = h(\mathbf{x}, \mathcal{N}):=\operatorname{argmin}\left\{\|\mathbf{y}-\mathbf{x}\|_{\infty}: 
\mathbf{y} \in \mathcal{N}, \mathbf{y}(2)>\mathbf{x}(2)\right\} .
\end{align}
In other words, $h(\bfx)$ is the closest Poisson point w.r.t. the $\ell_\infty$ norm with strictly higher $y$ coordinate.
Also,  $h(\bfx)$ is well-defined for all $\bfx \in \bbR^2$.
If the underlying point process is clear from the context, then we would drop the second coordinate and simply denote it by $h(\mathbf{x})$. Note that $h(\mathbf{x})$ has been defined for all $\bfx \in \mathbb{R}^2$. The (two-dimensional) $\ell_\infty$ DSF is defined as the random graph $\mathcal{T}$ with vertex set $\mathcal{N}$ and edge set $\mathcal{E} := \{ \langle \mathbf{x} , h(\mathbf{x}) \rangle : \mathbf{x} \in \mathcal{N} \}$. By construction, the $\ell_\infty$ DSF is a directed outdegree-one graph without any cycle or loop. 

Set $h^{0}(\mathbf{x}) = \mathbf{x}$ and for $k \geq 1$, let $h^{k}(\mathbf{x}) = h(h^{k-1}(\mathbf{x}))$ denote the $k$-th $\ell_\infty$ DSF step of $\mathbf{x}$. By joining successive DSF steps $h^k(\mathbf{x}) : k \geq 0$ by line segments, we obtain a continuous path $\pi^{\mathbf{x}}$ starting at $\mathbf{x}$. For $\nu, \delta >0$, integer $n \ge 1$ and
 a path $\pi^{\mathbf{x}}$ with starting time $\sigma_{\pi^{\mathbf{x}}}$, consider the diffusively scaled path $\pi^{\mathbf{x}}_n( \nu, \delta):[\frac{\sigma_{\pi^{\mathbf{x}}}}{n^2 \nu}, \infty] \to (-\infty, \infty)$ given by 
 \begin{equation}\label{def:scaled_path}
    \pi^{\mathbf{x}}_n(t) = \pi^{\mathbf{x}}_n( \nu, \delta)(t) := \frac{\pi^{\mathbf{x}}_n( n^2 \nu t)}{n \delta}.
 \end{equation}
 It is clear that the starting time of $\pi^{\mathbf{x}}_n$ is $\frac{\sigma_{\pi^{\mathbf{x}}}}{n^2 \nu}$. We have the following result:
 \begin{theorem}\label{thm:conv_bm}
     There exist $\nu = \nu( \lambda) >0$ and $\delta= \delta ( \lambda) >0$ such that as $n \to \infty$, $\pi^{\mathbf{0}}_n( \nu, \delta) $ converges in distribution to the standard Brownian motion starting from $\mathbf{0}$. 
 \end{theorem}
 For $\mathbf{x}, \mathbf{y} \in \mathbb{R}^2$ with $\mathbf{x}(2) = \mathbf{y}(2)$ we consider the DSF paths $\pi^{\mathbf{x}}, \pi^{\mathbf{y}}$ and their coalescing time is denoted by 
\begin{align}
\label{def:CoalesTime}
T^{\mathbf{x}, \mathbf{y}} := \inf\{ t > 0 : \pi^{\mathbf{x}}(\mathbf{x}(2) + t) = \pi^{\mathbf{y}}(\mathbf{y}(2) + t)\}.
\end{align}
The other main result of this paper obtains an optimal tail decay estimate for $T^{\mathbf{x}, \mathbf{y}}$ and thereby implies that the (two-dimensional) $\ell_\infty$ DSF is connected a.s.
\begin{theorem}
\label{thm:DSF_CoalTime_Decay}
There exists $C_0 > 0$, depending only on $\lambda > 0$, such that for all $t > 0$ we have
$$
\mathbb{P}(T^{\mathbf{x}, \mathbf{y}} > t) \leq  \frac{C_0 (| \mathbf{x}(1) - \mathbf{y}(1) | \vee 1)}{\sqrt{t}}. 
$$
\end{theorem}
Clearly, Theorem \ref{thm:DSF_CoalTime_Decay} implies that the coalescing time $T^{\mathbf{x}, \mathbf{y}}$ is finite a.s. Theorem \ref{thm:DSF_CoalTime_Decay} gives us that the $\ell_\infty$ DSF is connected a.s.

\begin{theorem}
\label{prop:DSF_Tree}
The two-dimensional $\ell_\infty$ DSF is connected a.s. 
\end{theorem}

We should mention here that the two-dimensional  $\ell_2$ DSF was introduced by Baccelli and Bordenave in \cite{MR2292589} and it was used as a tool to analyze the asymptotic properties of the radial spanning tree. Additionally, Baccelli et al. conjectured that the two-dimensional  $\ell_2$ DSF is connected a.s. This conjecture was later proved by Coupier and Tran in \cite{MR2999212} using a Burton-Keane type argument \cite{MR990777}. Baccelli et al. \cite{MR2292589} showed that under diffusive scaling, an $\ell_2$ DSF trajectory converges in distribution to a Brownian motion and they conjectured that entire scaled DSF network should converge to the Brownian web (BW), which can be described intuitively as a collection of independent coalescing Brownian motions starting from everywhere in $\mathbb{R}^2$. In a seminal work, Fontes et. al. \cite{MR2094432} characterized the BW as a random variable taking values in a Polish space and provided conditions for studying convergence to the BW. \cite{MR3644280} provides an excellent survey on BW. In a recent paper, Coupier et. al. \cite{MR4203342} proved the second conjecture.  

The tail decay estimate of coalescing time of two $\ell_2$ DSF paths is one of the key ingredients to prove convergence to the Brownian web and it was obtained in Theorem 5.1 in \cite{MR4203342}. Very recently, in an independent work \cite{garcia2025directed} Sanchez obtained this estimate for the two-dimensional $\ell_p$ DSF for all $ p \in [1 , \infty]$. We should mention here that in \cite{garcia2025directed}, Sanchez
developed a unifying argument to analyze the higher-dimensional $\ell_p$ DSF for dimensions $d \geq 3$ and for $ p \in \{1, 2, \infty\}$. This allowed him to establish a tree-forest dichotomy behavior for the $\ell_p$ DSF depending on dimensions.

In order to prove Theorem \ref{thm:DSF_CoalTime_Decay} for the $\ell_\infty$ DSF, we broadly follow the footsteps of \cite{MR4203342}. However, there are significant challenges to overcome. For the analysis of the $\ell_2$ DSF
in \cite{MR4203342}, one of the main building blocks is the construction of infinitely many renewal steps. Towards that, Lemma 3.2 of \cite{MR4203342} was a geometric result which was crucially used. It says that for the $\ell_2$ DSF, if we push the current moving vertex up to some height, precisely half of the height of the current history region, then there exists a cone of deterministic angular width, which avoids the history set. This observation was used in \cite{MR4203342} to bound the growth of history sets and to construct a good event, which would ensure a decrease in the height of history regions by a fixed amount and whose probability is uniformly bounded away from zero.    

For any $p \in [1,\infty)$
we consider  $\ell_{p}$ DSF, and observe that  by pushing the moving vertex above by a fraction $c_p \in (0, 1)$ (where $c_p$ depends only on $p$) of the height of the current history region, we obtain an unexplored cone of fixed angular width, and the same argument can be repeated. However,  for the $\ell_{\infty}$ DSF, no such $c_p < 1$ exists and this argument fails. 
As a result, we need new ideas to deal with this.

Before ending this section, we make the following remark. Since, Theorem \ref{thm:DSF_CoalTime_Decay} gives us the required estimate on the tail decay of coalescing time of two $\ell_\infty$ DSF paths. It is possible to follow the same methodology as in \cite{MR4203342} and show that the diffusively scaled $\ell_\infty$ DSF converges in distribution to the Brownian web. However, in this paper, we will not do that. We should mention here that in \cite{garcia2025directed}, Sanchez showed that the diffusively scaled $\ell_p$ DSF converges to BW for any $p \in [1, \infty]$. 

This paper is organized as follows. In Section \ref{section_2}, we introduce the joint exploration process that describes the movement algorithm governing the evolution of $k \geq 1$ DSF paths. In Section \ref{section_3}, we construct a sequence of renewal steps for the stochastic process associated with a single DSF path and analyze its properties at these renewal steps; this section also contains the proof of   Theorem~\ref{thm:conv_bm}. In Section \ref{sec:joint_renewal}, we define the renewal steps for the joint exploration process of two DSF paths, study the properties of the joint process at these renewal steps  and conclude the proofs of Theorem \ref{thm:DSF_CoalTime_Decay} and Theorem~\ref{prop:DSF_Tree}.

\section{Joint exploration process of DSF paths}
\label{section_2}

Fix $k\in \mathbb{N}$ and choose $k$ points $\mathbf{x}^{1}, \ldots, \mathbf{x}^{k}$ in $\mathbb{R}^2$. We will define the joint exploration process of DSF paths starting from $\mathbf{x}^{1}, \ldots, \mathbf{x}^{k}$ as $\{ (g_n(\mathbf{x}^1), \ldots, g_n(\mathbf{x}^k), H_n): n \ge 0\}$ which tracks the $k$ DSF paths in tandem. In what follows, for $\mathbf{x} \in \mathbb{R}^2$ we write $\|\mathbf{x}\|$ for $\|\mathbf{x}\|_{\infty}$ 
to simplify the notation. 

We need to introduce some notation.
For $r \in \mathbb{R}$ the (open) upper and lower half-planes are respectively defined as 
$$
\mathbb{H}^{+}(r) := \{ \mathbf{y} \in \mathbb{R}^2 : 
\mathbf{y}(2) > r\} \text{ and }
\mathbb{H}^{-}(r) := \{\mathbf{y} \in \mathbb{R}^2 : \mathbf{y}(2) < r\}.
$$
For $\mathbf{x} \in \mathbb{R}^2$ and $r \geq 0$ let 
$B (\mathbf{x},r ) := \{\mathbf{y} \in \mathbb{R}^2 : \|\mathbf{y} - \mathbf{x}\| < r \}$ 
denote the $\ell_\infty$ open ball of radius $r$ centered at $\mathbf{x}$ and 
$B^{+} (\mathbf{x},r ):= B (\mathbf{x}, r ) \cap \mathbb{H}^{+}(\mathbf{x}(2))$ denote the upper part of it.
 We take the starting configuration as $\mathbf{x}^{1}(2) = \cdots = \mathbf{x}^{k}(2)$ and $H_0 = \emptyset$. In general, the $k$ starting points can be general $k$ deterministic points in $\mathbb{R}^2$ or they can be Poisson points too. The starting set $H_0$ can be a non-empty set of the form $(\cup_{j=1}^{l} B^+(\mathbf{y}^j, r_j))\cap \mathbb{H}^{+}\bigl ( \min\{\mathbf{x}^i(2) : 1 \le i \le k \} \bigr )$ for some $\mathbf{y}^1, \cdots, \mathbf{y}^l \in \mathbb{R}^2$ and $r_1 , \cdots , r_l > 0$. We only require that the initial configuration $(\mathbf{x}^1, \cdots, \mathbf{x}^k, H_0)$ satisfies the following two conditions:
\begin{itemize}
    \item[(i)] $\mathbf{y}^j(2) < \min\{ \mathbf{x}^i(2) : 1 \le i \le k \} $ for all $1 \le j \le l$ and 
    \item[(ii)] $\mathbf{x}^i \notin H_0$
    for all $1 \le i \le k$.
\end{itemize}

W.l.o.g. we assume  $\mathbf{x}^{1}(2) = \cdots = \mathbf{x}^{k}(2) = 0$. We define the joint exploration process below in an inductive manner.  
For all $i=1,\ldots, k$ set $g_0(\mathbf{x}^i)= \mathbf{x}^i$ . 
To ensure that the paths move in tandem, only the point(s) with the least coordinate $y$, moves and all other points remain unchanged. 
In case several vertices have the least $y$ coordinate, they all move at once. For $1 \leq i \leq k$ we define $ g_1(\mathbf{x}^i )  :=  h(\mathbf{x}^i)$ .
Given $(g_1(\mathbf{x}^1), \cdots , g_1(\mathbf{x}^k))$,
set $r_1 := \min\{g_1(\mathbf{x}^i)(2) : 1 \leq i \leq k \}$ and we define the set $H_1 = H_1(\mathbf{x}^1, \cdots, \mathbf{x}^k)$ as  
$$
H_1 :=  \bigl ( \cup_{i=1}^k B^{+} ( \mathbf{x}^i,\| \mathbf{x}^i - h(\mathbf{x}^i)\| ) \bigr ) \bigcap \mathbb{H}^+(r_1).
$$
Given the first step $(g_1(\mathbf{x}^1), \cdots , g_1(\mathbf{x}^k))$, we define the `move' set after the first step as 
\begin{align*}
W^{\text{move}}_1 & := \{ g_1(\mathbf{x}^i) : g_1(\mathbf{x}^i)(2) = r_1 , 1\le i \le k\} .
\end{align*}
By the definition, we see that $W^{\text{move}}_1$ is a singleton set a.s.

More generally, for $n \geq 1$ given 
$(g_n(\mathbf{x}^1), \cdots , g_n(\mathbf{x}^k), H_n)$, we set $r_n := 
\min\{ g_n(\mathbf{x}^i)(2) : 1 \leq i \leq k \}$  and 
$W^{\text{move}}_{n} := \{ g_n(\mathbf{x}^i) : g_n(\mathbf{x}^i)(2) = r_n, 1\le i \le k\} $.
We define 
\begin{align*}
 g_{n+1}(\mathbf{x}^i ) & := 
 \begin{cases}
  h(g_{n}(\mathbf{x}^i )) & \text{ if } g_{n}(\mathbf{x}^i ) \in W_{n}^{\text {move }} \\ 
  g_{n}(\mathbf{x}^i ) & \text{ otherwise}. 
 \end{cases}
 \end{align*}
Note that $W^{\text{move}}_{n}$ is a singleton set a.s. and with a slight abuse of notation, it will be used to denote 
the moving vertex as well. Let $r_{n+1} := \min\{g_{n+1}(\mathbf{x}^i)(2) : 1 \leq i\leq k \}$ and the history set $H_{n+1}$ is defined as the region 
$$
H_{n+1} :=  \bigl ( H_n \cup B^{+} ( W^{\text{move}}_{n},\| W^{\text{move}}_{n} - h(W^{\text{move}}_{n}) \| ) \bigr ) 
\bigcap \mathbb{H}^+(r_{n+1}).
$$
This completes the description of the $(n+1)$-th step of the joint exploration process given as $(g_{n+1}(\mathbf{x}^1), \cdots , g_{n+1}(\mathbf{x}^k), H_{n+1}) $.

For any $n \geq 1$, given $(g_n(\mathbf{x}^1), \ldots, g_n(\mathbf{x}^k), H_n)$, to obtain the next step, the PPP in the unexplored part of the upper half-plane, i.e., in $\mathbb{H}^+(W^{\text{move}}_n(2)) \setminus H_n$ can be replaced with an independent PPP.  
As observed in \cite{MR4203342}, we can start with an i.i.d. collection $\{\mathcal{N}_n:n \in \mathbb{N}\}$ independent of the original PPP $\mathcal{N}$. For each $n \geq 1$ after the $n$-th step, in the unexplored part of the upper half-plane the point process can be replaced with $\mathcal{N}_{n+1}$ and it would not change the distribution of the next step. This gives us the so-called auxiliary exploration process, which has the same distribution as the joint exploration process. For more details, we refer the reader to \cite{MR4203342}. In what follows, we will actually work with the auxiliary exploration process, although we denote it using the same notation $\{(g_n(\mathbf{x}^1), \ldots, g_n(\mathbf{x}^k), H_n) : n \geq 0\}$. 

 Let 
\begin{align}
\label{def:Filtration_1}
\{ {\cF}_n = {\cF}_n(\mathbf{x}^1, \cdots, \mathbf{x}^k) := \sigma ((g_l(\mathbf{x}^1), \ldots, g_l(\mathbf{x}^k)) : 0\le l \le n) : n \geq 0 \}
\end{align}
denote the natural or minimal filtration with respect to which the exploration process is adapted.
The same argument as in Proposition 2.2. in \cite{MR4203342} proves the following proposition.
\begin{proposition}
  The process $\{ (g_n(\mathbf{x}^1), \ldots, g_n(\mathbf{x}^k), H_n) : n \geq 0\}$  is a $\{ {\cF}_n : n \geq 0 \}$ adapted Markov chain.  
\end{proposition}

\begin{figure}
    \centering
    \begin{tikzpicture}
        \draw (-0.8,0) rectangle (1.2,1);
         \draw (-0.3, 0.7) rectangle ( 2.7, 2.2 );
         \draw (-1.8, 1.8) rectangle (1.2, 3.3);
          \filldraw[fill = lightgray, draw= black] (-1.8, 2.6) rectangle ( 1.2, 3.3);
         \draw ( 0.5, 2.6) rectangle ( 1.9, 3.3);
        \draw [black,fill] (0.2,0) circle [radius=0.1] node[below = 2] {$\bfx$};
        \draw (0.2,0) -- ( 1.2, 0.7);
        \draw [black,fill] (1.2, 0.7) circle [radius=0.1];
        \draw [black, dashed] (1.2, 0.7) -- (3, 0.7) node[right = 2] {$g_1(\bfx)$};
         \draw  ( 1.2, 0.7) -- ( -0.3, 1.8);
        \draw [black,fill] ( -0.3, 1.8) circle [radius=0.1] node[below = 2] {$g_2(\bfx)$};
        \draw ( -0.3, 1.8) -- ( 1.2, 2.6);
         \draw [black,fill] ( 1.2, 2.6) circle [radius=0.1];
          \draw  ( 1.2,2.6) -- (0.7,3.3);
          \draw [black,fill] ( 0.7, 3.3) circle [radius=0.1] node[above = 2] {$g_4(\bfx)$};
          \draw [black, dashed] ( 1.2,2.6) -- (2.2,2.6) node[right = 2] {$g_3(\bfx)$};
    \end{tikzpicture}
    \caption{This picture illustrates the first $4$ steps of the marginal process $\{g_n( \bfx) : n \geq 0\}$. Grey region denotes the history set $H_3$. At the fourth step, we have $H_4 = \emptyset$ and a renewal event occurs as described in Section~\ref{section_3}.}
    \label{fig:renewal_step_single_path}
\end{figure}

\section{Renewal sequence for marginal path}
\label{section_3}

In this section, we define the sequence of  renewal steps for $k =1$, i.e., for the marginal process of a single DSF path $\{ g_n(\mathbf{x}^1) : n \geq 0\}$.  
Set $\beta^1_0 = 0$ and for $\ell \geq 1$ the $\ell$-th renewal step denoted by $\beta^1_\ell$ defined as 
\begin{align}
\label{def:Beta_Step}
\beta^1_\ell := \inf\{ n > \beta^1_{\ell - 1} : H_n = \emptyset \}.
\end{align}
For any $\ell \ge 1$, the r.v. $\beta^1_\ell$ is a stopping time w.r.t. the filtration $\{{\cF}_{n} : n \geq 0\}$.
Note that the above definition of marginal renewal step makes sense only for the $\ell_\infty$ DSF as it is impossible to have such a step for the $\ell_p$ DSF with $p \in [1,\infty)$. In fact, even for the $\ell_\infty$ DSF, for the joint process of $k =2$ paths, we cannot have such a step. For $k \geq 2$, the joint renewal step will be defined differently (see Section \ref{section_3}). 

Proposition \ref{prop:Beta_Tail_1} implies that $\beta^1_\ell < \infty$ a.s. for all $\ell \geq 1$. 

\begin{proposition}
\label{prop:Beta_Tail_1}
There exist $c_0, c_1 > 0$ depending only on $\lambda > 0$ such that  for all $\ell \geq 0$  and $n \in \bbN$ we have 
\begin{align}
 \label{eq:tau_Tail}
 \mathbb{P}(\beta^1_{\ell + 1} - \beta^1_{\ell} > n \mid {\cF}_{\beta^1_{\ell}}) \leq c_0 \exp{(-c_1 n)}
.
\end{align}
\end{proposition}

We comment that in Proposition \ref{prop:Beta_Tail_1} and in several other places in this paper we deal with universal decay constants. This means that their values would depend only on the PPP intensity $\lambda$ and on $k$, the number of DSF paths considered. Typically, they are denoted by $c_0, c_1 $, but their values may change from one line to another. Proposition \ref{prop:Beta_Tail_1} will be proved through a sequence of lemmas. Before that we need to introduce some notations.

For any bounded subset $H \subset \mathbb{R}^2$ we define it's height as $L(H):=\sup \{\mathbf{y}(2)-\mathbf{x}(2): \mathbf{x}, \mathbf{y} \in H\}$. 
We set $L(\emptyset )=0$. Given the $n$-th step $(g_n(\mathbf{x}^1), H_n)$, the newly generated history rectangle is defined as 
$$
H^{\text{new}}_{n+1} := B^+(g_n(\mathbf{x}^1), \|g_n(\mathbf{x}^1) -  g_{n+1}(\mathbf{x}^1)\|).
$$
To show that the history set becomes empty, it is equivalent to show that the process $\{ L_n := L(H_n) : n \geq 0\}$ hits zero. The next corollary is straightforward to observe and its proof has been left to the reader.

\begin{corollary}
\label{cor:History_Height}  
For each $n \geq 1$ almost surely we have
\begin{itemize}
    \item[(i)] $H_{n+1} =  (H^{\text{new}}_{n+1} \cup H_n ) \cap \mathbb{H}^+(g_{n+1}(\mathbf{x}^1)(2))$ and 
    \item[(ii)] 
    $L_{n + 1} < L_n \vee L(H^{\text{new}}_{n+1})$.
\end{itemize}
\end{corollary}

For any set $A \subset \mathbb{R}^2$, let $\rho (A)$ denote the boundary of the set. For $n \geq 1$ let $\rho^+_{n+1}$ denote a subset of the boundary set $\rho(H_{n+1}^{\text{new}})$ defined as 
$$
\rho^+_{n+1} := \rho(H_{n+1}^{\text{new}}) \cap \mathbb{H}^+(g_{n}(\mathbf{x}^1)(2)).
$$
Recall that $\mathbb{H}^+(g_{n+1}(\mathbf{x}^1)(2))$ denotes the open upper half-plane.  The `top' part of this boundary set $\rho^+_{n+1}$ is defined as 
$$
\rho^+_{T, n+1} :=  \rho^+_{n+1} 
\cap  \{ \mathbf{y} \in \mathbb{R}^2 :  \mathbf{y} (2) = 
g_{n}(\mathbf{x}^1)(2) + \|g_{n}(\mathbf{x}^1) - g_{n+1}(\mathbf{x}^1)\| \}.
$$
We refer the reader to Figure \ref{figure1}.

\begin{figure} 
\begin{center}
\begin{tikzpicture}
\draw [lightgray,fill] (0,0) rectangle (3,1.5);
\draw [red, ultra thick] (3,0) -- (3,1.5);
\draw [red, ultra thick] (0,0) -- (0,1.5);
\draw [black] (1.5,0) -- (2.4,1.5);
\draw [blue, ultra thick] (0,1.5) -- (3,1.5);
\draw [black, thick] (-1,0) -- (4,0);
\draw [black,fill] (1.5,0) circle [radius=0.1] node[below=2] {$g_n(\bfx^1)$} ;
\draw [black,fill] (2.4,1.5) circle [radius=0.1] node[above=2] {$g_{n + 1}(\bfx^1)$};

\end{tikzpicture}
\caption{This figure illustrates the occurrence of up step. The blue part denotes $\rho^+_{T, n+1}$ and union of blue part and red part denote  $\rho^+_{ n+1}$. The grey region represents the region $H^\text{new}_{n+1}$. }
\label{figure1}
\end{center}
\end{figure}

We now introduce two special types of steps, referred to as the `top' step and the `up'
step, which will play a key role in proving Proposition~\ref{prop:Beta_Tail_1}.
\begin{definition}[`Top' step and `up' step for $k=1$]
\label{def:top_step}
Starting from $\mathbf{u} \in \mathbb{R}^2$, we say that the DSF step $h(\mathbf{u})= h(\mathbf{u}, \mathcal{N})$ is an up step if we have $ h(\mathbf{u})(2) = \mathbf{u}(2) + \|h(\mathbf{u}) -\mathbf{u}\|$.

As a continuation of this, we say that the $(n+1)$-th step is an `up' step if it belongs to the top part of the boundary, viz., $
g_{n+1}(\mathbf{x}^1) \in \rho^+_{T, n+1}$.

The $(n+1)$-th step is called a `top'  step if it is an up step and we also have $\|g_{n+1}(\mathbf{x}^1) - g_{n}(\mathbf{x}^1)\| \geq 1$.
\end{definition}

 In addition, if this step belongs to the top part of the boundary set, then it is called an up step. The inequality in \eqref{eq:TopStep_HistHeight} explains the significance of top step. For any event $A$ (belongs to some appropriate $\sigma$-algebra), we denote 
 the indicator function as $\mathbf{1}_A$.
For any $n \geq 0$ we have 
\begin{align}
\label{eq:TopStep_HistHeight}  
\mathbf{1}_{\{(n+1)\text{-th step is a top step}\}}L_{n+1} \leq (L_n - 1 )\vee 0 \text{  a.s},
\end{align}
i.e., whenever top step occurs the height of the history set must decrease by at least one.
We need to show that for any $n \geq 1$ given $\mathcal{F}_n$, the probability that the $(n+1)$-th step is a top step is uniformly bounded from below. Lemma \ref{lem:UpStepBound} helps us to achieve this objective. 
  
\begin{lemma}
    \label{lem:UpStepBound}
We have the following:
\begin{itemize}
    \item[(i)] $\mathbb{P}(h(\mathbf{0}, \mathcal{N}) \text{ is up step}) = 1/2$.

    \item[(ii)] Consider $\mathbf{v}^1, \cdots ,\mathbf{v}^l \in \mathbb{H}^-(0)$  and $r_1, \cdots, r_l > 0$ such that 
    $\mathbf{0} \notin \cup_{j = 1}^{l} B^+(\mathbf{v}^j, r_j)$ for all $1 \leq j \leq l$. For any $r > 0$ we have 
    \begin{align*}
    & \mathbb{P} \Bigl( h(\mathbf{0}, \mathcal{N}) \text{ is up step} \mid \bigl( \cup_{j=1}^l B^+(\mathbf{v}^j, r_j) \cup B^+(\mathbf{0}, r) \bigr ) \cap \mathcal{N} = \emptyset \Bigr ) \geq 1/2.  
    \end{align*}
   \item[(iii)] For any $r_1 \in (r, 3r/2)$ we have 
    \begin{align*}
    & \mathbb{P} \Bigl( h(\mathbf{0}, \mathcal{N}) \text{ is up step} \mid  \bigl( \cup_{j=1}^l B^+(\mathbf{v}^j, r_j) \cup B^+(\mathbf{0}, r) \bigr ) \cap \mathcal{N} = \emptyset,  B^+((0,r), r_1 - r)
    \cap \mathcal{N} \neq  \emptyset  \Bigr ) \\
    & \geq 1/2.  
    \end{align*}     
\end{itemize}
\end{lemma}

\begin{proof} For any subset $A \subset \mathbb{R}^2$, let $T(A)\subset \mathbb{R}^2$ defined  as $T(A) := \{(\mathbf{x}(2),\mathbf{x}(1)) : \mathbf{x} \in  A \}$ denotes reflection about the line $y = x$. Clearly, $\mathcal{N}^\prime := T(\mathcal{N})$ has the same distribution  as $\mathcal{N}$ (see \cite[Chapter 5]{MR3791470}). 
For Item (i), we observe that for any realization $\omega$ of $\mathcal{N}$ belonging to the event $\{ h(\mathbf{0}, \mathcal{N}) \text{ is up step} \}$, the realization $T(\omega)$ in $\mathcal{N}^\prime$ belongs to the complementary event $\{ h(\mathbf{0}, \mathcal{N}^{\prime}) \text{ is
not an up step} \}$. On the other hand, for $\omega^\prime$  a realization of $\mathcal{N}^{\prime}$ belonging to the event $\{ h(\mathbf{0}, \mathcal{N}^{\prime}) \text{ is
not an up step} \}$, the realization $T^{-1}(\omega)$ in 
$\mathcal{N}$ belongs to the event $\{ h(\mathbf{0}, \mathcal{N}) \text{ is up step} \}$. Therefore, we have       
 $$
 \mathbb{P}(h(\mathbf{0}, \mathcal{N}) \text{ is up step}) =  \mathbb{P}(h(\mathbf{0}, \mathcal{N}^\prime) \text{ is not an up step}) = \frac{1}{2}.
 $$
 
For Item (ii), on the event $\{B^+(\mathbf{0}, r)\cap \mathcal{N} = \emptyset \}$, we must have $\|h(\mathbf{0})\| > r$. Further, the distribution of the PPP $\mathcal{N}$ in the region $(\mathbb{H}^+(0) \setminus B^+(\mathbf{0}, r))$ is independent of the event $ \{B^+(\mathbf{0}, r)\cap \mathcal{N} = \emptyset \}$. We consider  
any realization $\omega^\prime$ in $\mathcal{N}^\prime$ that belongs to the event 
$$
\{ h(\mathbf{0}, \mathcal{N}^\prime
\setminus (\cup_{j=1}^l B^+(\mathbf{v}^j, r_j) \cup B^+(\mathbf{0}, r))) \text{ is not an up step}\}.
$$
As $\mathbf{v}^j(2) < 0$ for all $ 1 \leq j \leq l$ with $ \mathbf{0} \notin \cup_{j=1}^{l} B^+(\mathbf{v}^j, r_j)$, the geometry of these $\ell_\infty$ balls ensures that in the  transformed realization $T^{-1}(\omega^\prime)$ belongs to the event  
$$
\{h \bigl ( \mathbf{0}, \mathcal{N} \setminus (\cup_{j=1}^l B^+(\mathbf{v}^j, r_j) \cup B^+(\mathbf{0}, r)) \bigr ) \text{ is an up step } \}.
$$
However, the converse does not necessarily hold. Therefore, we have
\begin{align*}
  &  \bbP \Bigl( h(\mathbf{0}, \mathcal{N}) \text{ is up step} \mid \bigl( \cup_{j=1}^l B^+(\mathbf{v}^j, r_j) \cup B^+(\mathbf{0}, r) \bigr ) \cap \mathcal{N} = \emptyset \Bigr ) \\
  & \ge \bbP \Bigl( h(\mathbf{0}, \mathcal{N}^\prime) \text{ is not an up step} \mid \bigl( \cup_{j=1}^l B^+(\mathbf{v}^j, r_j) \cup B^+(\mathbf{0}, r) \bigr ) \cap \mathcal{N}^\prime = \emptyset \Bigr ).
\end{align*}
This completes the proof for Item (ii).

The same argument as in Item (ii) applies for Item (iii) as well.  This completes the proof.
\end{proof}

\begin{corollary}
\label{cor:TopStepBound}
There exists $p_1 \in (0,1)$, which depends only on $\lambda$, such that for any $n \geq 0$  we have 
$$
\mathbb{P} (L_{n+1} \leq  (L_n - 1)\vee 0 \mid {\cF}_n) \geq p_1.
$$
\end{corollary}
\begin{proof}
Using the inequality in \eqref{eq:TopStep_HistHeight}  we obtain
\begin{align}
\label{eq:TopStepBound}
& \mathbb{P} (L_{n+1} \leq  (L_n - 1)\vee 0 \mid {\cF}_n) \nonumber \\
& \geq \mathbb{P} ((n+1)\text{-th  step is a top step} \mid {\cF}_n) \nonumber \\
& = \mathbb{P} ((n+1)\text{-th step is a top step} \mid g_n(\mathbf{x}^1), H_n) \nonumber  \\
& = \mathbb{P} ((n+1)\text{-th  step  is an up step} 
\mid B^+(g_n(\mathbf{x}^1), 1)\cap \mathcal{N}_{n+1}  = \emptyset, g_n(\mathbf{x}^1), H_n) \nonumber \\
& \qquad \qquad  \times \mathbb{P}(B^+(g_n(\mathbf{x}^1), 1)\cap \mathcal{N}_{n+1}  = \emptyset) \nonumber \\
& \geq \frac{1}{2} \exp{(-2\lambda)}.
\end{align}
The step in \eqref{eq:TopStepBound} follows from  Item (ii) of Lemma \ref{lem:UpStepBound}. This completes the proof.
\end{proof}

Next, we need to bound the increase of the process $\{L_n : n \geq 0\}$. For any $n \geq 1$ let $g_n^\uparrow := g_n(\mathbf{x}^1) + (0,L_n)$. For $\bfx \in \bbR^2$ and $A \subset \bbR^2$, we define $\bfx \circ A := \{ \bfx + \bfy : \bfy \in A \} \subset \bbR^2$. Let $R_{n+1}$ denote the r.v. defined as
\begin{align}
    \label{def:R_n}
    R_{n+1} := \inf \bigl \{ l \geq 0 : \bigl( g_n^\uparrow \circ \bigl( [-l, l] \times [0,1] \bigr ) \bigr ) \cap \mathcal{N}_{n+1} \neq \emptyset \bigr \}.
\end{align}
In other words, the r.v. $R_{n+1}$ represents the smallest width $l$ so that the rectangle $g_n^\uparrow \circ ([-l, l] \times [0,1])$ contains $\mathcal{N}_{n+1}$ point(s). The rectangle $g_n^\uparrow \circ ([- L_n - 1,  L_n  + 1] \times [0,1])$ is shown in Figure \ref{figure2}. We observe that for all $l > 0$ we have $ g_n^\uparrow \circ \bigl( [-l, l] \times [0,1] \bigr ) \cap H_n = \emptyset$. We observe that the distribution of $R_{n+1}$ does not depend on $\mathcal{F}_{n}$.
\begin{lemma}
\label{lem:ubound_Lprocess}
For any $n \geq 0$ and for any $m \geq 1$ we have 
\begin{align*}
    \mathbb{P}(L_{n+1} > L_{n} + m \mid \mathcal{F}_{n}) \leq \exp{(-2 \lambda (L_{n} + m))}.
\end{align*}
\end{lemma}

\begin{proof} Using Corollary \ref{cor:History_Height} we obtain 
\begin{align}
\label{eq:ubound_Lprocess}
    \mathbb{P}(L_{n+1} > L_{n} + m \mid \mathcal{F}_{n})     & \le \mathbb{P}(L(H^{\text{new}}_{n+1}) > L_{n} + m \mid \mathcal{F}_{n}) \le \mathbb{P}( R_{n+1} \ge L_{n} + m ) .
\end{align}
For the last step in \eqref{eq:ubound_Lprocess} we note 
that for any $m \geq 1$, the region
$g^\uparrow_n \circ \bigl( [-L_n - m, L_n + m ] \times [0, 1]\bigr )$ avoids $H_n$ and it is contained in the ball $B^+(g_n(\mathbf{x}^1), L_n + m)$. Therefore, to occur the event
$\{ L(H^{\text{new}}_{n+1}) > L_n + m \}$, the rectangular region $g^\uparrow_n \circ \bigl( [-L_n - m, L_n + m ] \times [0, 1]\bigr )$ must be free of points from $\mathcal{N}_{n+1}$. As a result, \eqref{eq:ubound_Lprocess} follows. The proof follows from the distribution of $R_{n+1}$ which is independent of $\mathcal{F}_{n}$. 
\end{proof}

Lemma \ref{lem:ubound_Lprocess} controls the probability $\mathbb{P}(L_{n+1} > L_{n} + m \mid \mathcal{F}_{n})$, for 
any $m \ge 1$. The next lemma controls the probability when the amount of increase is less than one, i.e., $\mathbb{P}(L_{n+1} \in ( L_{n}, L_{n} + 1) \mid \mathcal{F}_{n})$.
It is important to observe that for
$ L_{n+1} \in (L_n, L_n +1) $, it is not necessary to have $R_{n+1} > L_n$ and the inequality in \eqref{eq:ubound_Lprocess_1} does not always hold.
We have described such a situation in Figure \ref{figure2}.

\begin{lemma}
\label{lem:ubound_Lprocess_1}
Fix any $n \geq 1$. There exists $m_0 \ge 2$ such that given $L_n \ge m_0$, we have 
\begin{align}
\label{eq:ubound_Lprocess_1}
    \mathbb{P}(L_{n+1} \in (L_{n}, L_{n} + 1) \mid \mathcal{F}_{n}) \le \frac{2}{ L_n }.
\end{align}
\end{lemma}

\begin{figure}
\begin{center}
\begin{tikzpicture}
\draw [lightgray,fill] (0,0) rectangle (4,3);
\draw [lightgray,fill] (0,0) rectangle (-4,3);
\draw [black, thick] (0,0) -- (4,0);
\draw [black, thick] (0,0) -- (-4,0);
\draw [black, thick] (4,1) -- (4,3) -- (3,3);
\draw [black, thick] (-4,1) -- (-4,3) -- (-3,3);
\draw [black, dashed] (0,0) -- (0,3);
\draw (-5,3) rectangle (5,4);
\draw[black, thick] (-4, 4) -- (-4, 3);
\draw[black, thick] (-3, 4) -- (-3, 3);
\draw[black, thick] (-2, 4) -- (-2, 3);
\draw[black, thick] (-1, 4) -- (-1, 3);
\draw[black, thick] (-0, 4) -- (-0, 3);
\draw[black, thick] (1, 4) -- (1, 3);
\draw[black, thick] (2, 4) -- (2, 3);
\draw[black, thick] (3, 4) -- (3, 3);
\draw[black, thick] (4, 4) -- (4, 3);
\draw (4,0) -- (6,0);
\draw (-4,0) -- (-6,0);
\draw (4,0) rectangle (5,1);
\draw (-4,0) rectangle (-5,1);
\draw [black,fill] (0,0) circle [radius=0.1] node[below=2] {$g_n(\bfx^1)$} ;
\draw [black,fill] (0,3) circle [radius=0.1] ;
\draw [black,fill] (4.5,0.7) circle [radius=0.1] ;
\draw [black,fill] (-4.3,0.2) circle [radius=0.1] ;
\draw [black,fill] (2,3.8) circle [radius=0.1] ;
\draw[arrow] (2.8, 3.6) -- (3, 4.3) node[above] {$g_n^\uparrow \circ ([ - L_n - 1, L_n + 1] \times [0,1])$};
\draw[arrow] (4.6, 0.4) -- (5.5, 0.5) node[right] {$\Box^r_n $};
\draw[arrow] (-4.6, 0.4) -- (-5.5, 0.5) node[left] {$\Box^l_n $};
\draw [black, dashed] (-4.3,0) -- (-4.3, 3.3) -- (4.3, 3.3) -- (4.3,0);
\draw[arrow] (1,0) -- (1,3) node[midway, right]{$L_n$};
\draw[arrow] (0,3) -- (0,4.2) node[above]{$g_n^\uparrow $};
\draw ( 0,0) -- (-4.3, 0.2);
\end{tikzpicture}
\caption{The point $g_n(\bfx^1)$ may connect with a Poisson point in the region $\Box^l_n \cup \Box^r_n$  rather than connecting with a Poisson point in the region $g_n^\uparrow \circ ([ - L_n - 1, L_n + 1] \times [0,1])$. As a result, the event $E_{ n + 1}$ occurs, and the height of the history set may increase but the amount of increase is bounded by 1.}
\label{figure2}
\end{center}
\end{figure}

\begin{proof}
We first observe that
\begin{align}
\label{eq:L_nIncrease_Bd_1}
 & \mathbb{P}(L_{n+1} \in (L_{n}, L_{n} + 1) \mid \mathcal{F}_{n})   \nonumber\\
  & \leq  \mathbb{P} \Bigl ( L_{n+1} \in (L_{n}, L_{n} + 1) , R_{n+1} \le L_n + 1 \mid \mathcal{F}_{n} \Bigr ) +
 \mathbb{P} \bigl ( R_{n+1} > L_n + 1\mid \mathcal{F}_{n}  ).
\end{align}
The second term in  \eqref{eq:L_nIncrease_Bd_1} can be bounded by $\exp{(-2 \lambda (L_n + 1))}$.
To deal with the first term in \eqref{eq:L_nIncrease_Bd_1}, we consider two square boxes of unit length $\Box^r_n$ and $\Box^l_n$  (see Figure \ref{figure2}) defined respectively as
\begin{align*}
 \Box^r_n & := (g_n(\mathbf{x}^1) + (L_n, 0)) \circ ([1,0]\times [0,1]) \text{ and }   \\ 
 \Box^l_n & := (g_n(\mathbf{x}^1) + (-L_n, 0)) \circ ([-1,0]\times [0,1]).
\end{align*}
Considering the above two boxes we define the event $E_{n+1}$ as  
\begin{align}
\label{def:E_n}
E_{n+1} & := \bigl \{ h \bigl( g_n(\mathbf{x}^1), \mathcal{N}_{n+1} \setminus B^+(g_n(\mathbf{x}^1), L_n) \bigr ) \in \Box^r_n\cup \Box^l_n \bigr \} \text{ and } \nonumber \\
F_{n+1} & := \bigl \{ h \bigl ( g_n(\mathbf{x}^1), \mathcal{N}_{n+1} \setminus B^+(g_n(\mathbf{x}^1), L_n) \bigr ) \in B^+(g_n(\mathbf{x}^1), L_n + 1) \bigr \}. 
\end{align}
We observe that on the event $\{ R_{n+1} \le L_n + 1 \}$, in order to have $L_{n+1} \in (L_{n}, L_{n} + 1) $, the event $E_{n+1} \cap F_{n+1}$ necessarily occurs. The occurrence of the event $E_{n+1} \cap F_{n+1}$ may not be enough to ensure an increase in $L_{n+1}$, as the boxes
$\Box^r_n$ and $\Box^l_n$ may still be covered (at least partially) by the history set $H_{n}$. In fact, connection to a $\mathcal{N}_{n+1}$ point in $\Box^r_n \cup \Box^l_n$ may not result an increase in $L_{n+1}$. Therefore, we have
$$
\mathbb{P} \Bigl ( L_{n+1} \in (L_{n}, L_{n} + 1) , R_{n+1} \le L_n + 1 \mid \mathcal{F}_{n} \Bigr ) \le \mathbb{P}(E_{n+1}\cap F_{n+1}\mid \mathcal{F}_{n}  )
\le \mathbb{P}(E_{n+1} \mid  F_{n+1}, \mathcal{F}_{n}  ).
$$

In order to bound the probability $\mathbb{P}(E_{n+1} \mid  F_{n+1}, \mathcal{F}_{n} )$,
we observe that the PPP $\mathcal{N}_{n+1}$ is independent of $\mathcal{F}_{n}$ and henceforth is distributed independently in the region $B^+(g_n(\mathbf{x}^1), L_n + 1) \setminus B^+(g_n(\mathbf{x}^1), L_n)$.
We consider square boxes $\Box^\uparrow_{n, j} := g_n^\uparrow \circ 
([j, j+1 ]\times [0,1] )$ for all $j \in [- L_n - 1 , L_n] \cap \mathbb{Z} $ (see Figure \ref{figure2}). Given that the event $F_{n+1}$ has occurred, all the boxes $\Box^r_{n}, \Box^l_{n}$ and $\Box^\uparrow_{n, j}$ for $j \in [L_n - 1 ,  L_n ]\cap \mathbb{Z}$, are equally likely to hold the nearest $\mathcal{N}_{n+1}$ point (w.r.t. the $\ell_\infty$ norm). Therefore, \eqref{eq:L_nIncrease_Bd_1} becomes  
\begin{align*}
\mathbb{P} \Bigl ( L_{(n+1)} \in (L_{n}, L_{n} + 1) \mid \mathcal{F}_{n} \Bigr )  
 & \leq \mathbb{P} ( E_{n+1} \mid F_{n+1}, \mathcal{F}_{n+1} ) + \mathbb{P} ( R_{n+1} > L_n + 1 )    \\
& \leq \frac{1}{L_n }  + \exp{(-2 \lambda(L_n + 1 ))} \leq \frac{2}{ L_n  },
\end{align*}
for all large $L_n$. This completes the proof.
\end{proof}

We are now ready to prove Proposition \ref{prop:Beta_Tail_1}. For $M \ge m_0$ where $m_0$ is as in Lemma \ref{lem:ubound_Lprocess_1},  we define $\tau^M := \inf \{ n \ge 1 : L_n \le M \}$. Clearly, $\tau^M$ is a stopping time w.r.t. the filtration $\{ \mathcal{F}_n : n \geq 0\}$.

\begin{lemma}
\label{lem:exp_tail_tau_1}
    There exist constants $c_0 , c_1 > 0$ depending only on $\lambda$ and $M$ such that for any $n \in \bbN $ we have 
    \[
    \bbP( \tau^M > n ) \le c_0 \exp{( - c_1 n)}.
    \]
\end{lemma}
Before proving Lemma \ref{lem:exp_tail_tau_1}, we first show that Lemma \ref{lem:exp_tail_tau_1} gives us Proposition \ref{prop:Beta_Tail_1}.
\vspace{0.5cm}

\noindent {\bf Proof of Proposition~\ref{prop:Beta_Tail_1}:}
Using Item (ii) of Lemma \ref{lem:UpStepBound} we obtain 
\begin{align*}
& \mathbb{P}( L_{n+1} = 0 \mid L_n \leq M, \mathcal{F}_n) \\
& \ge \mathbb{P}\bigl ( (n+1)\text{-th step is top step } \mid B^+(g_n(\mathbf{x}^1), M) \cap \mathcal{N}_{n+1} = \emptyset, L_n \leq M, \mathcal{F}_n \bigr ) \times \\
& \qquad \qquad \mathbb{P}\bigl ( B^+(g_n(\mathbf{x}^1), M) \cap \mathcal{N}_{n+1} = \emptyset \mid L_n \leq M, \mathcal{F}_n \bigr )\\
& \geq \frac{1}{2}\exp{(-2\lambda M^2)}.    
\end{align*}
Therefore, Proposition \ref{prop:Beta_Tail_1} follows from Lemma \ref{lem:exp_tail_tau_1}.
\qed
\vspace{0.5cm}

\noindent {\bf Proof of Lemma \ref{lem:exp_tail_tau_1}:} We first show that there exist $ \alpha, \epsilon > 0$ and $M \ge m_0 > 1$ such that the process 
$$
\{ S_{n\wedge \tau^M} := \exp{(\alpha L_{n \wedge \tau^M} + \epsilon (n\wedge \tau^M))} : n \geq 0\}
$$
gives an $\{ \mathcal{F}_n : n \geq 0\}$ adapted non-negative supermartingale. 
For $\tau^M \leq n$, it is straightforward to see that $\mathbb{E}[ S_{(n+1)\wedge \tau^M} \mid \mathcal{F}_n] = S_{n \wedge \tau^M}$ a.s. On the other hand, for $\tau^M > n$, 
we choose $\alpha >0$ small enough so that $\mathbb{E}[ \exp{(\alpha R_{n+1})}] < \infty$ and we obtain 
\begin{align}
\label{eq:mar_ineq_1}
& \mathbb{E}[ \exp{(\alpha(L_{n+1} - L_{n}))} \mid \mathcal{F}_n] \nonumber\\
&= \mathbb{E}[ \exp{(\alpha(L_{n+1} - L_{n}))} (\mathbf{1}_{\{L_{n+1} > L_n + 1\}} + \mathbf{1}_{\{ L_{n+1} \in (L_n, L_n + 1)\}} + \mathbf{1}_{\{ L_{n+1} \in (L_n - 1, L_n)\}} \nonumber \\
&\quad \quad+ \mathbf{1}_{\{ L_{n+1} \leq L_n - 1\}}) \mid \mathcal{F}_n] \nonumber \\
& \leq \mathbb{E}[ \exp{(\alpha R_{n+1})\mathbf{1}_{\{R_{n+1} > M + 1\}}} \mid \mathcal{F}_n] + 
\exp{(\alpha)} \mathbb{P}( L_{n+1} \in (L_n, L_n + 1) \mid \mathcal{F}_n) +\nonumber \\
& \qquad \qquad \qquad  \bbP( L_{n+1} \in (L_n - 1, L_n) \mid \cF_n)+ \exp{(-\alpha)} \mathbb{P}( L_{n+1} \leq (L_n - 1) \mid \mathcal{F}_n) \nonumber \\
& \leq \mathbb{E}[ \exp{(\alpha R_{n+1})\mathbf{1}_{\{R_{n+1} > M + 1\}}} \mid \mathcal{F}_n] + 
\exp{(\alpha)} \mathbb{P}( L_{n+1} \in (L_n, L_n + 1) \mid \mathcal{F}_n) + \nonumber \\
& \qquad  \qquad  \bbP( L_{n+1} \in (L_n - 1, L_n) \mid \cF_n)+ \exp{(-\alpha)} \mathbb{P}( (n+1)\text{-th step is top step} \mid \mathcal{F}_n) \nonumber \\
& \leq \mathbb{E}[ \exp{(\alpha R_{n+1})\mathbf{1}_{\{R_{n+1} > M + 1\}}}] + 
 \frac{2 \exp{(\alpha)}}{M} + \big( 1 - \frac{2}{M} - \frac{\exp(-2 \lambda)}{2} \big)+ \frac{\exp{(-(\alpha + 2 \lambda))}}{2} \nonumber \\
&= \mathbb{E}[ \exp{(\alpha R_{n+1})\mathbf{1}_{\{R_{n+1} > M + 1\}}}] + 
 \frac{2 \exp{(\alpha)}}{M} + \big( 1 - \frac{2}{M} - \frac{\exp(-2 \lambda)}{2} ( 1 - \exp(-\alpha)) \big),
\end{align}
where the last inequality in \eqref{eq:mar_ineq_1} follows from Lemma \ref{lem:ubound_Lprocess_1} and Corollary \ref{cor:TopStepBound}.
It is not difficult to see that 
we can choose $M \geq m_0$ large such that \eqref{eq:mar_ineq_1} becomes strictly lesser than  $1$. That would mean  that the process $\{ S_{n \wedge \tau^M} : n \geq 0\}$ (for the chosen $\alpha > 0$ and $M \ge m_0$) is a non-negative supermartingale. Hence,  there exists a r.v. $Z$ with $\mathbb{E}(Z) < \infty$ such that as $n \to \infty$, we have
$$
S_{n \wedge \tau^M} \to Z \text{ a.s.}
$$
Assuming $S_0 > M$ we obtain 
\[
\mathbb{E} \big[ \mathbf{1}_{ \{\tau^M = \infty \}} \exp{(\epsilon n)} \big] \leq \mathbb{E} \big[ \mathbf{1}_{ \{ \tau^M = \infty \}} S_n \big] \leq \mathbb{E} \big[ S_{n \wedge \tau^M } \big] \leq M \text{ for all }n\ge 1.
\]
This gives us $\mathbb{P}( \tau^M = \infty) \leq M \exp{(-\epsilon n)}$ for all $n \ge 1$, which implies $\mathbb{P}(\tau^M  < \infty) = 1$. 
Therefore, the limiting r.v. $Z$ must be $
Z = S_{\tau^M}$ a.s. Hence, we obtain  $$
\mathbb{E}[\exp{(\epsilon \tau^M )}] \leq \mathbb{E}[S_{\tau^M }] \leq S_0 .
$$
Finally, Markov inequality gives us
\[
\mathbb{P} (\tau^M  > n) = \mathbb{P} ( \exp{(\epsilon \tau^M )} > \exp{(\epsilon n)} ) \leq \exp{(-\epsilon n)} \mathbb{E} [ \exp{(\epsilon \tau^M )}] \leq 
M \exp{(-\epsilon n)} .
\]
\qed

Now consider the process $\{ g_{\beta^1_{\ell}} (\mathbf{x}^1): \ell \ge 0 \}$, where the renewal steps $\beta^1_{\ell}$ have 
been defined in \eqref{def:Beta_Step}. Let us introduce the 
following random variable
\begin{equation}\label{def:Y_ell}
    Y_{\ell} = Y_{\ell}(\mathbf{x}^1) := g_{\beta^1_{\ell}}  (\mathbf{x}^1) \ \ \text{ for } \ell \ge 0.
\end{equation}
Since at the $\beta^1_{\ell}$-th step the history set becomes empty and the translation invariant nature of our model ensure that:
\begin{proposition}\label{prop:iid_incre}
    The process $\{ Y_{\ell} - Y_{\ell -1 } : \ell \ge 2 \}$
    is a sequence of i.i.d. random vectors whose distribution does not depend on the starting point $\mathbf{x}^1$. 
Moreover, the distribution of $(Y_{2} - Y_{1})(1) $ is symmetric about $0$. 
\end{proposition}

Next we show that the random variable $Y_{\ell}(1) - Y_{\ell - 1}(1) $ has moments of all
orders. Fix $\ell \ge 0$ and consider the random variable
\begin{equation}\label{def:W_ell}
    W_{\ell +1} = W_{\ell +1}(\mathbf{x}^1) := \sum_{m = \beta^1_{\ell}}^{\beta^1_{\ell +1} -1} \|g_m(\mathbf{x}^1) - h(g_m(\mathbf{x}^1)) \|.
\end{equation}

\begin{lemma}\label{lem:exp_tail_width}
  There exists constants $c_0,c_1>0$ depending only on $\lambda$ and $M$ such that for any $n \in \mathbb{N}$ and $\ell \ge 1$ we have
  \begin{equation}
  \label{eq:exp_tail_width}
    \bbP(\|(Y_{\ell + 1} - Y_{\ell})\| > n \mid \cF_{\beta^1_{\ell}}) \le \bbP( W_{\ell+1} >n \mid \cF_{\beta^1_{\ell}}) \le c_0 \exp(-c_1 n^{1/2}).
  \end{equation}
 
\end{lemma}

\begin{proof}
The first inequality in \eqref{eq:exp_tail_width}
is straightforward to observe. For the next part, we  work 
conditionally on $\cF_{\beta^1_{\ell}}$. We recall the definition
of $ R_{n +1}$  (see \eqref{def:R_n}). We claim  that given $\cF_{\beta^1_{\ell}}$ for any $m \ge 0$, we have
\begin{equation}\label{eq:claim}
      \max_{0 \le n \le m } \bigl(\|g_{\beta^1_{\ell} + n}(\mathbf{x}^1) - h(g_{\beta^1_{\ell} +n}(\mathbf{x}^1)) \| \bigr) + 1 \le \max_{0 \le n \le m} (R_{\beta^1_{\ell} + n + 1} + 1).
\end{equation}
We will prove \eqref{eq:claim} by induction. We first consider $m = 0$. On the event $\{\|h(g_{\beta^1_{\ell}}(\mathbf{x}^1)) - g_{\beta^1_{\ell}}(\mathbf{x}^1)\|\le 1\}$, clearly \eqref{eq:claim} holds. On the other hand, while dealing with  the event $\{\|h(g_{\beta^1_{\ell}}(\mathbf{x}^1)) - g_{\beta^1_{\ell}}(\mathbf{x}^1)\| >  1\}$, we note that $L_{\beta^1_{\ell}} = 0$. Therefore, the argument of Lemma \ref{lem:ubound_Lprocess} gives us that 
$$ 
\|h(g_{\beta^1_{\ell}}(\mathbf{x}^1)) - g_{\beta^1_{\ell}}(\mathbf{x}^1)\| + 1 \le  R_{\beta^1_{\ell} + 1} + 1  .
$$
This proves \eqref{eq:claim} for the choice $m = 0$. 

 We assume that \eqref{eq:claim} holds for some $m \ge 0$. Given that $g_{\beta^1_{\ell}}(\mathbf{x}^1) = \mathbf{w}$, we first consider the event
$\{ \|h^{m+1}(\mathbf{w}) - h^{m}(\mathbf{w})\| \le 
L_{\beta^1_{\ell} + m} + 1 \}$. By induction hypothesis, we have 
\[
L_{\beta^1_{\ell} + m}   + 1 \le 
\max_{0 \le n < m} \|h^{n+1}(\mathbf{w}) - h^n(\mathbf{w})\| + 1
\le
\max_{0 \le n \le m} (R_{n + 1} + 1) 
\]
and therefore, on the event
$\{ \|h^{m+1}(\mathbf{w}) - h^{m}(\mathbf{w})\| \le 
L_{\beta^1_{\ell} + m} + 1 \}$
\eqref{eq:claim} holds for $m+1$. On the other hand, on the event 
$\{ \|h^{m+1}(\mathbf{w}) - h^{m}(\mathbf{w})\| > L_{\beta^{1}_\ell + m} +  1\}$, the argument of Lemma \ref{lem:ubound_Lprocess} ensures that we must have 
\begin{align*}
\|h^{m+1}(\mathbf{w}) - h^{m}(\mathbf{w})\| \le R_{\beta^1_{\ell} + 1} \text{ implying }\|h^{m+1}(\mathbf{w}) - h^{m}(\mathbf{w})\| + 1 \le R_{\beta^1_{\ell} + 1} + 1 .  
\end{align*}
As a result, \eqref{eq:claim} holds for $m + 1$. This proves \eqref{eq:claim} by the induction method.  

Now \eqref{eq:claim} gives us that 
\begin{align}\label{eq:claim_1}
    & \bbP( W_{\ell+1} >n \mid \cF_{\beta^1_{\ell}}) \le \bbP \Big( \sum_{m = 1}^{\beta^1_{\ell+1} - \beta^1_\ell} (R_{\beta^1_\ell + m} + 1)> n \mid \cF_{\beta^1_{\ell}} \Big) \nonumber \\
        & \le \bbP \Big( \sum_{m = 0}^{\lfloor n^{1/2}\rfloor}  \max_{0 \le j \le m} R_{j + 1} > n -  n^{1/2} \Big) + \bbP ( \beta^1_{\ell + 1} - \beta^1_{\ell + 1} > n^{ 1/2} \mid \cF_{\beta^1_{\ell}}).
\end{align}
Thanks to Proposition~\ref{prop:Beta_Tail_1}, the second term of the left-hand side of 
\eqref{eq:claim_1} is bounded above by $c_0 \exp(-c_1 n^{1/2})$. We consider the first term:
\begin{align*}
     \bbP \Big( \sum_{m = 0}^{\lfloor n^{1/2}\rfloor}  \max_{0 \le j \le m} R_{j + 1} > n - n^{1/2}\Big) & \le \bbP \Big( \sum_{m = 0}^{\lfloor n^{1/2}\rfloor}  \max_{0 \le j \le m} R_{j + 1} > n/2 \Big)\\
     & \le \bbP \Big( \bigcup_{m =0}^{\lfloor n^{1/2}\rfloor} R_{m+1} > c n^{1/2}  \Big)  \\
     &\le (\lfloor n^{1/2}\rfloor+1) \bbP( R_1 >  cn^{1/2}  )  \\ 
     &\le (\lfloor n^{1/2}\rfloor+1) c_0 \exp{(- c_1 n^{1/2})}.
\end{align*}
The last step follows from the distribution of $R_1$. 
\end{proof}

\noindent{\bf Proof of Theorem~\ref{thm:conv_bm}:}
   The time for the $\ell$-th renewal, captured by the distance 
   traveled by the process in second co-ordinate, is given by
   \begin{align}\label{def:ren_time}
       T_{\ell} = T_{\ell}(\mathbf{x}^1) := g_{\beta^1_{\ell}}(\mathbf{x}^1)(2) - \mathbf{x}^1(2) = Y_{\ell}(2) -  Y_{\ell - 1}(2) \text{ for any } \ell \ge 1.
   \end{align}
Observe that $T_{\ell + 1} - T_{\ell} \le W_{\ell +1}$ for all $\ell \ge 1$ and hence, all moments of $ T_{\ell + 1} - T_{\ell} $ exist.

Now for the path $\pi^{\mathbf{0}}$ we define the modified path $\widetilde{\pi}$ which is linear between successive 
renewal steps of $\pi^{\mathbf{0}}$. Proposition~\ref{prop:iid_incre} and Lemma~\ref{lem:exp_tail_width} allow us to deduce that the displacements between successive renewal steps are i.i.d. with mean zero and exponential moments. Therefore, we can apply Donsker's invariance principle (see \cite[Page 90, Section 8]{MR1700749}) to conclude that the diffusively scaled  modified path converges to the standard Brownian motion for the choice $\nu = \bbE(T_2 - T_1)$ and $\delta = \sqrt{\text{Var}(Y_{2}(1) - Y_1(1))}$.
Finally, we need to show that the modified path and the actual path has same diffusive scaling limit. Mathematically, it is enough to show that for any $s>0$ and $\varepsilon >0$,
\begin{align}\label{eq:modi_path}
    & \bbP ( \sup \{ | \widetilde \pi_n (t) - \pi^{\mathbf{0}}_n(t) | : t \in [0,s] \} > \varepsilon ) \nonumber \\
    &= \bbP ( \sup \{ | \widetilde \pi (t) - \pi^{\mathbf{0}}(t) | : t \in [0,n^2 \nu s] \} > n \delta \varepsilon ) \to 0 \text{ as } n \to \infty, 
\end{align}
where $\widetilde \pi_n $ denotes the modified scaled path. 

We recall the definition of $W_{\ell +1}(\mathbf{0})$ from \eqref{def:W_ell}. Note
that for any $t \in [T_{\ell}, T_{\ell +1}]$ and for  all $\ell \ge 1$ we have $|\widetilde{\pi}(t) - \pi^{\mathbf{0}}(t)| \le W_{\ell +1}$. Since we can have at most $\lfloor n^2 \nu s \rfloor$ many renewals in time $n^2 \nu s$, we have
\begin{align*}
    &\bbP ( \sup \{ | \widetilde \pi (t) - \pi^{\mathbf{0}}(t) | : t \in [0,n^2 \nu s] \} > n \delta \varepsilon ) \\
     & \le \bbP( \max \{W_{\ell +1} : \ell =0, \ldots, \lfloor n^2 \nu s \rfloor \} > n \delta \varepsilon ) \\
   & \le (\lfloor n^2 \nu s \rfloor +1) c_0\exp{(-c_1 (n \delta \varepsilon)^{1/2})} \to 0 \text{ as } n \to \infty,
\end{align*}
where the last inequality follows from Lemma~\ref{lem:exp_tail_width}. This completes the proof. 
\qed


\begin{figure}
    \centering
    \begin{tikzpicture}
    \tikzstyle{arrow_beidseitig}=[thick,<->,>=stealth]
    \filldraw[lightgray, draw = black] (-1.2, 0.7) rectangle ( 3.4,3);
    \filldraw[lightgray, draw = black] (5.6, 0.7) rectangle (10, 3);;
        \filldraw[gray, draw= black] ( -1.5, 0.7) rectangle ( -0.9,1);
        \filldraw[gray, draw= black] ( -0.3, 0.7) rectangle ( 1.1,1.4);
        \filldraw[gray, draw= black] ( 7.3, 0.7) rectangle ( 8.6,1.2);
        \begin{scope}[fill opacity=1.2]
          \fill[white] (0.1,3) -- (1.1, 2) -- (2.1,3) -- cycle; 
        \end{scope}
        \begin{scope}[fill opacity=1.2]
          \fill[white] (6.9,3) -- (7.9, 2) -- (8.9,3) -- cycle; 
        \end{scope}
         \draw [black,fill] (1.1,0.7) circle [radius=0.1] node[below = 2]{$g_{\beta_ 1}(\bfx^1)$};
          \draw[black,fill] (7.9,1.2) circle[radius=0.1] node[ below =2] {$g_{\beta_1}(\bfx^2)$} ;
         \draw [black, thick] ( -2, 0.7) -- (11,0.7);
          \draw [black, thick] ( -2, 2) -- (11, 2);
         \draw[red, fill] (1.1, 2) circle[radius = 0.1] node[below = 2] {$g^{\uparrow, 1}_{\tau_1}$};
        \draw[red, fill] (7.9, 2) circle[radius = 0.1] node[below = 2] {$g^{\uparrow, 2}_{\tau_1}$};
        \draw (-1.2, 0.7) rectangle ( 3.4, 3);
        \draw (5.6, 0.7) rectangle (10, 3);
        \draw (0.1, 2) rectangle ( 2.1, 3);
        \draw (6.9, 2) rectangle (8.9, 3);
        \draw (1.1,2) -- (0.1,3);
        \draw (1.1,2) -- (2.1,3);
        \draw (7.9,2) -- (8.9,3);
        \draw (7.9,2) -- (6.9, 3);
       \draw[arrow_beidseitig] (-1.6, 0.7) -- (-1.6, 3) node[midway, left]{$\kappa +1$};
       \draw[arrow_beidseitig] (11, 0.7) -- (11, 2) node[midway, left]{$\kappa$};
       \draw[arrow_beidseitig] (-1.2, 3.1) -- (3.4, 3.1) node[midway, above]{$B^+( g^{\downarrow, 1}_{\tau_1}, \kappa +1 )$};
        \draw[arrow_beidseitig] (5.6, 3.1) -- (10, 3.1) node[midway, above]{$B^+( g^{\downarrow, 2}_{\tau_1}, \kappa +1 )$};
        \draw[black, fill] (1, 2.6) circle[radius = 0.1];
        \draw[black, fill] (8, 2.6) circle[radius = 0.1];
        \draw (1, 2.6) -- (1.1, 0.7);
        \draw (7.9, 1.2) --  (8, 2.6);
        \draw (1, 2.6) --  (1.1, 2);
        \draw  (8, 2.6) -- (7.9, 2);
    \end{tikzpicture}
    \caption{This picture illustrates joint renewal step of the joint exploration process $\{g_n(\bfx^1), g_n(\bfx^2), H_n): n \ge 0 \}$. Red dots represent the projected vertices $g^{\uparrow, 1}_{\beta_1}$ and $g^{\uparrow, 2}_{\beta_1}$, i.e., positions of restart. The grey region denotes the history set $H_{\beta_1}( \bfx^1, \bfx^2)$.}
    \label{fig:enter-label}
\end{figure}

\section{Regeneration for the Joint process of \texorpdfstring{$k=2$}{k = 2} paths}
\label{sec:joint_renewal}

For $k \geq 2$ we consider the joint exploration process. It follows that for any $n \geq 1$, as long as the set $\{ g_n(\mathbf{x^1}) , \cdots , g_n(\mathbf{x^k})\} $ contains multiple points, the history set $H_n$ can never be empty. 
Therefore, for the joint exploration process of general $k\geq 2$ DSF paths, we require a modified definition of the renewal step. In the next subsection,
we first define a sequence of `good' steps such that we have a good control over the height of generated history regions.

\subsection{`Good' steps for joint process of \texorpdfstring{$k\ge2$}{k >= 2} paths}
\label{subsec:joint_good}

In order to define the sequence of good steps, for a technical reason (see the proof of Lemma \ref{lem:small_increase_G_n}) instead of working with the height process $\{L_n : n \geq 0\}$, we consider a related  process $\{G_n : n \geq 0\}$ which is defined as 
\[
G_n := L_n +  \# W^{\text{stay}}_n,
\]
where for any set $A$, we denote the cardinality by $\#A$. We consider this process in multiples of $k$ steps, i.e, $\{ G_{nk} : n \ge 0\}$. 
Set $\tau_0 = 0 $. For $j \geq 1$ we define 
$$
\tau_j := \inf \{ nk > \tau_{j-1} : G_{nk} \leq \kappa, W^{\text{move}}_{nk}(2) > W^{\text{move}}_{\tau_{j-1}}(2) + \kappa + 1 \},
$$
where $\kappa \ge (k-1) + m_0$ ($m_0$ is same as in Lemma~\ref{lem:ubound_Lprocess_1}) is any positive constant. In what follows, these steps will be referred as `good' steps as we have control over the heights of history regions on these steps.
For any $j \geq 1$ the r.v. $\tau_j$ is a stopping time w.r.t. the filtration $\{ {\cF}_n : n \geq 0\}$.   
The next proposition shows that the r.v. $\tau_j$ is finite a.s. for all $j \geq 1$. 
\begin{proposition}
    \label{prop:Tau_tail_k}
   There exist constants $c_0, c_1 > 0$  depending only on $\lambda$ and  $k$ such that for all $j \geq 1$ and $n \in \mathbb{N}$ we have
    \begin{align}
    \label{eq:Tau_tail_k}
    \mathbb{P}(\tau_{j+1} - \tau_j > n \mid {\cF}_{\tau_j}) \le c_0 \exp( - c_1 n). 
    \end{align}
\end{proposition}

To prove Proposition \ref{prop:Tau_tail_k}, we first modify the definition of `top' step and `up' step for the joint process of $k\geq 2$ DSF paths. The newer history rectangle created  due to 
the $(n+1)$-th step is given as $H^{\text{new}}_{n+1} := B^+(W^{\text{move}}_n, \|W^{\text{move}}_n - h(W^{\text{move}}_n)\|)$. Recall that for $n \geq 1$, the set $W^{\text{move}}_n$ is a singleton set a.s. and with a slight abuse of notation, it is used to denote the moving vertex as well. The analogous boundary sets $\rho^+_{n+1}$ and the $\rho^+_{T,n+1}$ are respectively defined as 
\begin{align*}
 \rho^+_{n+1} & := \rho(H^{\text{new}}_{n+1}) \cap \mathbb{H}^+(W^{\text{move}}_n(2)) \text{ and }\\  
 \rho^+_{T, n+1} & := \rho^+_{n+1} \cap \{ \mathbf{y} : \mathbf{y}(2) = 
 W^{\text{move}}_n(2) + \|W^{\text{move}}_n - h(W^{\text{move}}_n)\| \}.
\end{align*}
Top step and up step for $k \geq 2$ DSF paths are defined as follows.

\begin{definition}(`Top' step and `Up' step for $k \geq 2$) 
\label{def:Top_Step}
We say that the $(n+1)$-th step is an up step if 
$h(W^{\text{move}}_n) \in \rho^+_{T, n+1}$. 

Additionally, the $(n+1)$-th step is called a top  step if  one of the following conditions is satisfied: 
\begin{itemize}
    \item[(i)] either $(n+1)$-th step is an \textit{up} step and $$ h(W^{\text{move}}_n)(2) \in \big[ W^{\text{move}}_n (2) + 1 , W^{\text{move}}_n (2) + L_n + 1/2 \big] \text{ or} 
    $$
    \item[(ii)]  $h(W^{\text{move}}_n) \in \{ g_n(\mathbf{x}^i) : 1\leq i \leq k \}$ and $ \|W^{\text{move}}_n - h(W^{\text{move}}_n)\|  
 \le L_n + 1/2 $ .
\end{itemize}
\end{definition}

The next lemma highlights the utility of the ` top steps to control the decay of the process $\{ G_{nk} : n \geq 0\}$.

\begin{lemma}
 \label{lem:Ln_Decrease}
 Given $L_{nk} > \kappa$,
on  $k$ consecutive occurrences of top steps we have
 \begin{align}
     \label{eq:Ln_Decrease}
G_{(n+1)k}\mathbf{1}_{ \{ (nk+j)\text{-th step is top step for all }1 \leq j \leq k\}} \leq G_{nk} - 1/2 \text{ a.s.}
  \end{align}
\end{lemma}
\begin{proof}
We  actually prove \eqref{eq:Ln_Decrease} in terms of the $L_{nk}$ function, i.e., we will show that
 \begin{align}
     \label{eq:Ln_Decrease_1}
L_{(n+1)k}\mathbf{1}_{ \{ (nk+j)\text{-th step is top step for all }1 \leq j \leq k\}} \leq L_{nk} - 1/2 \text{ a.s.}
  \end{align}
  We observe that (\ref{eq:Ln_Decrease_1}) immediately gives us \eqref{eq:Ln_Decrease}.

   For simplicity, we write the proof for $k = 2$. The argument remains the same for general $k \ge 2$. In case $g_{2n}(\mathbf{x}^1) = g_{2n}(\mathbf{x}^2)$, then (\ref{eq:Ln_Decrease}) readily follows from our earlier observation on the marginal process of a single DSF path. Therefore, we assume $g_{2n}(\mathbf{x}^1) \neq  g_{2n}(\mathbf{x}^2)$. Considering the case of two consecutive occurrences of the top steps with $L(H^{\text{new}}_{2n+1}) \leq L(H_{2n}) = L_{2n} $, we observe that we should necessarily have
   \begin{align}
   \label{eq:L_Increase}    
    L(H^{\text{new}}_{2(n+1)})\vee L(H^{\text{new}}_{2n+1}) \leq L_{2n} + 1/2.
   \end{align}
  We see from \eqref{eq:L_Increase} that the `possible' increase is bounded by at most $1/2$, whereas two consecutive occurrences of the top steps ensure that 
  $W^{\text{move}}_{2(n+1)}(2) \geq W^{\text{move}}_{2n}(2) + 1$
  resulting in a decrease of at least $1$ in the height function of older history rectangles. Therefore, (\ref{eq:Ln_Decrease}) follows. 

  We now consider the sub case $L(H^{\text{new}}_{2n+1}) > L_{2n}$. In case we have $ L(H^{\text{new}}_{2(n+1)}) \leq L_{2n}$, then (\ref{eq:L_Increase}) still holds and therefore, the same argument proves (\ref{eq:Ln_Decrease}). 
  On the other hand, considering the subcase $L(H^{\text{new}}_{2n+1})\wedge L(H^{\text{new}}_{2(n+1)}) > L_{2n}$ the `possible' amount of increase is bounded by at most $1/2 + 1/2 = 1$, however, for this sub case we must have $W^{\text{move}}_{2(n+1)}(2) \geq W^{\text{move}}_{2n}(2) + L_{2n}$. 
  As a result, older history height function must decrease by $L_{2n} > m_0 \ge 2$. This proves  (\ref{eq:Ln_Decrease}).
\end{proof}

Next, we show that for every $k$-th step of the joint process and given any history, the probability of $k$ consecutive occurrences of top step is uniformly bounded away from zero.
\begin{lemma}
 \label{lem:top_Step_LowerBound}
 Given $L_{nk} \geq \kappa \ge 2$, there exists $p_k > 0 $ depending on $\lambda$ such that we have 
\begin{align}
\label{eq:JtOccurGoodEvent}    
\mathbb{P}( (nk + j)-\text{th  step is a top step for all }1 \le j \le k \mid \mathcal{F}_{nk} ) \geq 
 p_k.
\end{align}
\end{lemma}
\begin{proof}
    For $1 \le j \le k$ we define the following two events: 
    \begin{align*}
    A_{nk + j} & :=  \{ B^+(W^{\text{move}}_{ nk + j-1} + (0,L_{nk+ j-1}),\frac{1}{8}) \cap \cN_{nk +  j} \neq \emptyset\}\text{ and }\\  
    B_{ nk + j} & := \{  B^+(W^{\text{move}}_{ nk+ j-1 },1) \cap \cN_{nk + j} = \emptyset \}.
    \end{align*}
    We observe that for any $1 \leq j \leq k$ by definition, the events $A_{ nk + j}$ and  $B_{ nk + j}$ are independent.   
    Using these events for $j = 1$ we obtain
    \begin{align}
    \label{eq:lem:top_Step_LowerBound}
     & \mathbb{P}( (nk+1)\text{-th step is a top step} \mid \mathcal{F}_{nk} ) \nonumber\\
     &\ge \mathbb{P}( (nk+1)\text{-th step is a top step}, B_{nk+1} \cap A_{nk + 1} \mid \mathcal{F}_{nk} ) \nonumber\\
      & = \mathbb{P} (  (nk+1)\text{-th step is a top step} \mid B_{ nk + 1} \cap  A_{ nk + 1},  \mathcal{F}_{nk} )\mathbb{P} (   B_{ nk + 1} \cap  A_{ nk + 1} \mid \mathcal{F}_{nk} ) \nonumber \\
         & = \mathbb{P} (  (nk+1)\text{-th step is a top step} \mid B_{ nk + 1} \cap  A_{ nk + 1},  \mathcal{F}_{nk} )\mathbb{P}(  A_{ nk + 1}  \mid \mathcal{F}_{nk} )
         \mathbb{P} ( B_{ nk + 1}  \mid \mathcal{F}_{nk} ) \nonumber \\
    & \ge \frac{1}{2} \mathbb{P}(  A_{ nk + 1}  \mid \mathcal{F}_{nk} )
         \mathbb{P} ( B_{ nk + 1}  \mid \mathcal{F}_{nk} ) \nonumber \\,
         & = \frac{1}{2} \mathbb{P}(  B^+(\mathbf{0},\frac{1}{8})\cap \cN_{nk + 1} \neq \emptyset )
         \mathbb{P} ( B^+(\mathbf{0},1) \cap \cN_{nk + 1} = \emptyset ) .
         \end{align}
     The last inequality in (\ref{eq:lem:top_Step_LowerBound}) follows from Item (iii) of Lemma \ref{lem:UpStepBound} and (\ref{eq:lem:top_Step_LowerBound}) follows from the translation invariance nature of the PPP $\cN_{nk + 1}$. Therefore, the r.h.s  in (\ref{eq:lem:top_Step_LowerBound})  can be bounded by 
     $\frac{1}{2} \exp{(- 2\lambda )}(1- \exp{(- \frac{\lambda}{32})})$. In the above argument, we have also used that given $\mathcal{F}_{nk}$, the events $A_{nk+1}$ and $B_{nk+1}$ are independent of $\mathcal{F}_{nk}$. 
 The auxiliary construction of the joint exploration process allows us to use the above argument repeatedly for $k$ many times. Then setting $p_k := \frac{1}{2^k} \exp{(- 2k\lambda )}(1- \exp{(- \frac{\lambda}{32})})^k >0 $ completes the proof.
\end{proof}

Lemma \ref{lem:top_Step_LowerBound} and Lemma \ref{lem:Ln_Decrease} together give us the following corollary.
\begin{corollary}
 \label{cor:Gn_Decrease}
 For any $n \geq 1$ we have 
 $$
 \mathbb{P}( G_{(n+1)k} \leq G_{nk} - 1/2 \mid {\cF}_{nk} ) \geq p_k,
 $$
 where $p_k > 0$ is as in Lemma \ref{lem:top_Step_LowerBound}.
\end{corollary}

We proceed as in the previous section. We first need to control the increase of the process $\{ G_{nk} : n \geq 0\}$. Let $W^\uparrow_{n}$ denote the point $W^\uparrow_{n} := W^{\text{move}}_{n} + (0, L_{n})$.  With a slight abuse of notation let $R_{n+1}$ be a non-negative r.v. defined for the joint exploration process of $k\geq 2$ DSF paths as
\begin{align}
\label{def:R_n_k}
R_{n + 1} := \inf\{ l > 0 : (W^\uparrow_{n} \circ ([-l,l]\times [0,1])) \cap \mathcal{N}_{n + 1} \neq \emptyset\}.    
\end{align}
The auxiliary construction process ensures that 
$\{R_{n} : n \geq 1\}$ forms a collection of i.i.d. random variables with exponentially decaying tails. 

\begin{lemma}
\label{lem:G_n_large_Increase_bd}
Given $G_{nk} > \kappa $, there exists constants $c_0, c_1 >0$ depending only on $\lambda$ such that for any $n \geq 0$  we have 
\begin{align*}
& \mathbb{P}\bigl( (G_{(n+1)k} - G_{nk}) > k \mid \mathcal{F}_{nk}) = \mathbb{P}\bigl( (L_{(n+1)k} - L_{nk}) > k \mid \mathcal{F}_{nk}) \le k c_0\exp{(-c_1 L_{nk})}.
\end{align*}

\end{lemma}

\begin{proof} The first inequality is straightforward to observe. Hence, we obtain 
\begin{align}
\label{eq:L_increase_morethan1} 
\{ L_{(n+1)k} - L_{nk} > k \} \subseteq  \cup_{j = 1}^k \{   L(H^{\text{new}}_{nk + j}) \ge L_{nk} + 1  \}.
\end{align}
The observation in (\ref{eq:L_increase_morethan1}) together with the same argument as in Lemma \ref{lem:ubound_Lprocess} give us the following inclusion relation:
\begin{align}
\label{eq:R_increase_morethan1} 
    \{ L_{(n+1)k} - L_{nk} > k \} \subseteq  \cup_{j = 1}^k \{   R_{nk + j} \ge L_{nk} + 1  \}.
    \end{align}
    The auxiliary construction of the joint exploration process, ensures that the family $\{ R_{nk+j} : 1 \leq j \leq k\}$ forms an i.i.d. collection of r.v.'s independent of $\mathcal{F}_{nk}$. Therefore, the proof follows by applying  union bound in (\ref{eq:R_increase_morethan1}). 
\end{proof}

Finally, Lemma \ref{lem:small_increase_G_n} obtains an upper bound for the probability           that the increase in $G_{nk}$ is bounded by 
$k$.

\begin{lemma}
\label{lem:small_increase_G_n}
    Fix any $n \ge 0$. Given $G_{nk} > \kappa$, we have 
\begin{align}
\label{eq:small_increase_G_n}    
 \bbP( G_{( n + 1)k} - G_{nk} \in ( 0, k) \mid \cF_{nk}) \le k \frac{2}{ L_{nk} }.
\end{align}
\end{lemma}

\begin{proof}
The argument is very similar to that of 
Lemma \ref{lem:ubound_Lprocess_1} and we provide a brief sketch here. We require minor
modifications since, conditional to $\mathcal{F}_{nk}$, the locations of the
vertices $\{g_{nk}(\mathbf{x^i}) : 1 \le i \le k\}$ are deterministic and the move vertex $ W^{\text{move}}_{nk}$ may connect to other vertices in the stay set. This may lead to a bounded amount of increase in the $L_{nk}$ function which cannot be controlled or dominated through the corresponding $R_{nk + 1}$ random variable. However, on such a step due to coalescence, the cardinality of the stay set decreases by $1$ and therefore, this would not lead to an increase for the $G_{nk}$ function.  
This is precisely the reason why we need to consider the function $G_{nk}$ instead of the $L_{nk}$ function. 

We consider the case for $k =2$. The argument is the same for general $k \ge 2$. W.l.o.g. we assume that $g_{2n}(\mathbf{x}^1) \neq g_{2n}(\mathbf{x}^2)$. Otherwise, we need to deal with a single DSF path only and in such situation, applying Lemma \ref{lem:ubound_Lprocess_1} repeatedly gives us the required result.
We observe that 
\begin{align}
\label{eq:G_n_large_Increase_bd}
& \mathbb{P}( G_{2(n+1)} - G_{2n} \in (0,2)\mid \mathcal{F}_{2n})\nonumber\\
& \leq \mathbb{P} \bigl ( G_{2(n+1)} - G_{2n} \in (0,2) \cap \{ \cap_{j=1}^2 R_{2n + j} < L_{2n} + 1\}  \mid \mathcal{F}_{2n} \bigr ) + \mathbb{P}(  \cup_{j=1}^2 R_{2n + j} > L_{2n} + 1 ) .
\end{align}
From Lemma \ref{lem:G_n_large_Increase_bd} we get that the second term in \eqref{eq:G_n_large_Increase_bd} is upper bounded by $2 \exp{(-2\lambda(L_{2n} + 1))}$. For the first term in \eqref{eq:G_n_large_Increase_bd} we observe that 
for each $1 \le j \le 2$, the possible amount of increase for the $L_{2n+j}$ function at the $j$-th step is bounded by $1$. Therefore, $W^{\text{move}}_{2n + j}$ connects to a stay vertex, that does not result an increase in the $G_{2n+j}$ function. 

Given that $ g_{2n}(\mathbf{x}^1) \neq g_{2n}(\mathbf{x}^2)$, on the event 
$\cap_{j= 1}^2 \{R_{2n + j} \le L_{2n} + 1\} $ we modify the definitions of the square boxes in the proof of Lemma \ref{lem:ubound_Lprocess_1}
as below:  
\begin{align*}
\Box^r_{2n + j} & := (W^{\text{move}}_{2n + j} + ( L_{2n + j} , 0 ) ) \circ ([1,0] \times [0,1] ) \text{ and} \\
\Box^l_{2n + j} & := (W^{\text{move}}_{2n + j} + ( -L_{2n + j} , 0 )) \circ ([-1,0] \times [0,1]),
\end{align*}
where $j = 0,1$.
Similarly,  for $ j \in \{1,2\}$ we modify the definitions of the events $E_{2n+j}$ and $F_{2n+j}$ as below:
\begin{align*}
E_{2n+j} & := \bigl \{ h \bigl( W^{\text{move}}_{2n + j-1}, \mathcal{N}_{2n+j} \setminus B^+(W^{\text{move}}_{2n+ j-1}, L_{2n}) \bigr ) \in \Box^r_{2n + j-1} \cup \Box^l_{2n + j-1} \bigr \} \text{ and }\\
F_{2n+j} & := \bigl \{ h \bigl( W^{\text{move}}_{2n + j-1}, \mathcal{N}_{2n+j} \setminus B^+(W^{\text{move}}_{2n+ j-1}, L_{2n}) \bigr ) \in B^+(W^{\text{move}}_{2n+ j-1}, L_{2n} + 1) \bigr ) \bigr \}.
\end{align*}
The above discussions suggest that given $\mathcal{F}_{2n}$ on the event 
$\cap_{j= 1}^2 \{R_{2n + j} \le L_{2n} + 1\}$
in order to have $G_{2(n+1)} - G_{2n} \in (0,2)$, the event $E_{2n+j}\cap F_{2n+j}$ must occur for some $ j =1,2$. Therefore,  the same argument of Lemma  \ref{lem:ubound_Lprocess_1} together with an union bound give us the required bound. This concludes the proof.
\end{proof}

We are now ready to complete the proof of Proposition \ref{prop:Tau_tail_k}. 
\vspace{0.5cm}

\noindent {\bf Proof of Proposition \ref{prop:Tau_tail_k}:}
The argument is similar to that of Proposition \ref{prop:Beta_Tail_1} and we provide a brief sketch here. Let 
$$
\widetilde{\tau} := \inf\{ n \ge 1 : G_{nk} \le \kappa \}, 
$$
where $\kappa \geq m_0 + (k-1)$ and $m_0 \ge 2$ is as in Lemma~\ref{lem:ubound_Lprocess_1}.
We observe that $k$ consecutive occurrences of top steps causes vertical advancement of at least $1$, i.e., given $G_{nk} > \kappa$ almost surely we have 
 \begin{align*}
W^{\text{move}}_{(n+1)k}(2)\mathbf{1}_{ \{ (nk+j)\text{-th step is top step for all }1 \leq j \leq k\}} \geq W^{\text{move}}_{nk}(2) + 1.
\end{align*}
As a result, in geometric steps the joint process crosses a height of $\kappa +1$. 
Therefore, to prove Proposition \ref{prop:Tau_tail_k} for $j = 0$ it suffices to show that the process 
$$
\{ S^2_{n \wedge \widetilde{\tau}} := \exp{ (\alpha G_{k(n \wedge \widetilde{\tau})}} + 
 \epsilon (n \wedge \widetilde{\tau})) : n \geq 0\}
$$
is a non-negative supermartingale for some $\alpha, \epsilon > 0$ and $ \kappa \geq m_0 + (k-1)$.

As observed earlier, for $\widetilde{\tau} \leq n$, the argument is straightforward. For $\widetilde{\tau} > n$
we define $ \widetilde{R}_{n+1} := \max\{ R_{nk + j} : 1 \le j \le k\}$ and observe that 
\begin{align}
\label{eq:Ineq_1}
& \mathbb{E}[ \exp{(\alpha(G_{(n+1)k} - G_{nk}))} \mid \mathcal{F}_{nk}] \nonumber\\
&= \mathbb{E}[ \exp{(\alpha(G_{(n+1)k} - G_{nk}))} ( \mathbf{1}_{\{G_{(n+1)k} > G_{nk} + k\}} + \mathbf{1}_{\{ G_{(n+1)k} \in (G_{nk}, G_{nk} + k)\}} + \nonumber \\ & \quad \quad \mathbf{1}_{\{ G_{(n+1)k} \in (G_{nk} -1/2, G_{nk}) \}} + \mathbf{1}_{\{ G_{(n+1)k} \leq G_{nk} - 1/2\}})  \mid \mathcal{F}_{nk}] \nonumber \\
& \leq \mathbb{E}[ \exp{(\alpha \widetilde{R}_{n+1})\mathbf{1}_{\{\widetilde{R}_{n+1} > L_{nk} +1\}}} \mid \mathcal{F}_{nk}] + 
\exp{(\alpha k)} \mathbb{P}( G_{(n+1)k} \in (G_{nk}, G_{nk} + k) \mid \mathcal{F}_{nk}) +\nonumber \\
& \qquad \qquad \qquad \bbP( G_{(n+1)k} \in (G_{nk} -1/2, G_{nk})) + \exp{(-\alpha/2)} \mathbb{P}( G_{(n+1)k} \leq G_{nk} - 1/2 \mid \mathcal{F}_{nk}) \nonumber \\
& \leq \mathbb{E}[ \exp{(\alpha \widetilde{R}_{n+1})\mathbf{1}_{\{\widetilde{R}_{n+1} > \kappa \}}}] + 
\Big[ k  \frac{2 \exp{(\alpha k)}}{\kappa} + p_k \Big] + \big( 1 - \frac{2 k}{\kappa} - p_k \big).
\end{align}
The last inequality in \eqref{eq:Ineq_1} follows from Lemma \ref{lem:G_n_large_Increase_bd} and Corollary \ref{cor:Gn_Decrease}.
The rest of the argument is same as that of Lemma \ref{lem:exp_tail_tau_1}. This completes the proof of Proposition \ref{prop:Tau_tail_k} for $j = 0$.

Now, let us prove the result for $j = 1$ and the  argument for general $j\geq 1$ is similar. We work conditionally on $\mathcal{F}_{\tau_1}$ and
observe that the ``new starting
configuration" $(g_{\tau_1} (\mathbf{x}^1), \cdots , g_{\tau_1} (\mathbf{x}^k), H_{\tau_1})$ satisfies the assumptions mentioned in Section \ref{section_2}. Therefore from the $\tau_1$-th step onwards, we can apply the same strategy and 
this compltes the proof.
\qed

\subsection{Renewal steps for the joint process of \texorpdfstring{$k=2$}{k=2} DSF paths}
\label{subsec:joint_renewal}

In this subsection, we define the joint renewal step for
$k \geq 2$ DSF paths. We need to introduce some notation first. By definition, for a $\tau_j$-th step we have 
$$
|g_{\tau_j}(\mathbf{x}^{i_1})(2) - g_{\tau_j}(\mathbf{x}^{i_2})(2)| \leq \kappa \text{ for all }1 \leq i_1 , i_2 \leq k.
$$
Let $C_{\frac{\pi}{4}} ( \mathbf{0}) := \{ r e^{i \theta} : r > 0, \theta \in [ \frac{\pi}{4} , \frac{3\pi}{4} ] \} $ be the cone with apex $\mathbf{0}$ and making an angle $\frac{\pi}{4}$ with the vertical axis. For $\mathbf{x} \in \bbR^2$, let $C_{\frac{\pi}{4}} ( \mathbf{x}) := \bfx \circ C_{\frac{\pi}{4}} ( \mathbf{0}) = \{ \mathbf{x} + \mathbf{y} : \mathbf{y} \in C_{\frac{\pi}{4}} ( \mathbf{0}) \}  $ denote the translated cone. We observe that if $h(\mathbf{x}) \in C_{\frac{\pi}{4}} ( \mathbf{x})$ then it must be a top step.  

For $1 \leq i \leq k $ let 
$g^{\downarrow,i}_j$ and $g^{\uparrow,i}_j$
respectively denote the projections of the point $g_{\tau_j}(\mathbf{x}^{i})$
on the lines $y = W^{\text{move}}_{\tau_j}(2)$ and $y = W^{\text{move}}_{\tau_j}(2)  + \kappa$. 
Considering the $\tau_j$-th step, the renewal event $\mathrm{Ren}_j = \mathrm{Ren}_j(\mathbf{x}^1, \cdots, \mathbf{x}^k)$ is defined as
\begin{align}
\label{def:E_j}
\mathrm{Ren}_j := \cap_{i=1}^k\{ & \# (B^+(g_j^{\downarrow, i}, \kappa +1) \cap \mathcal{N}_{\tau_j +1}) = \# \bigl ( B^+(g_j^{\uparrow, i},1) \cap C_{\pi/4}(g_j^{\uparrow, i}) \cap \mathcal{N}_{\tau_j +1} \bigr ) = 1 \}.
\end{align}

Set $\gamma_0 = 1$. Let $\gamma_\ell$ denote the number of good steps required for the $\ell$-th occurrence of the joint renewal event and it is defined as
\begin{align}
\label{def:Gamma_Renewal_k}
\gamma_\ell := \inf\{ j > \gamma_{\ell - 1} : \text{The event }\mathrm{Ren}_j \text{ occurs at the }j\text{-th good step}\}.
\end{align}
Let $\beta_\ell := \tau_{\gamma_\ell}$ denote the total number of steps required for the $\ell$-th occurrence of (joint) renewal event.
We observe that for any $j\geq 1$, the event $\{\gamma_\ell \le 
 j \}$ is not measurable w.r.t. $\mathcal{F}_{\tau_j}$ and therefore we need to extend this filtration as follows. 
For any $j\geq 1$ we define the $\sigma$-algebra $\mathcal{S}_j$ as 
\[
\mathcal{S}_j := \sigma \bigl ( \mathcal{F}_{\tau_j} , \mathrm{Ren}_1, \cdots , \mathrm{Ren}_j  \bigr ),
\]
and observe that for any $\ell \geq 1$, the r.v. $\gamma_\ell$ is a stopping time w.r.t. the filtration $\{ \mathcal{S}_j : j \geq 1\}$. This allows us to define the filtration 
\begin{align}
\label{def:G_ell}
\{ \mathcal{G}_\ell := \mathcal{S}_{\gamma_\ell} : \ell \geq 1 \},
\end{align}
and we observe that for any $\ell \geq 1$, the r.v.s $\gamma_\ell, \beta_\ell, g_{\beta_\ell}(\mathbf{x}^1), \cdots , g_{\beta_\ell}(\mathbf{x}^k)$ are all $\mathcal{G}_\ell$ measurable.

 In this paper, we need joint renewal steps for $k=2$ DSF paths only and we will prove the  next proposition, which gives us that $\beta_\ell$ is finite for all $\ell \geq 1$, for $k=2$ only. Its proof follows from the similar argument of Proposition 4.2 of \cite{MR4203342} . Basically, the auxiliary construction and bounded height of history regions at good steps together give us that at every good step, the probability of occurrence of a renewal event is uniformly bounded from below. 

\begin{proposition}
    \label{prop:Beta_Tail_k}
    There exist constants $c_0, c_{1} > 0$ depending only on $\lambda$ and $\kappa$ such that for all $n \in \mathbb{N}$ and $\ell \geq 0$  we have
    \begin{equation}\label{eq:beta_l}
        \mathbb{P}(\beta_{\ell +1 } - \beta_{\ell} > n \mid {\cG}_\ell) \leq c_0 \exp{(-c_{1} n)}.
    \end{equation}

\end{proposition}

To prove Proposition~\ref{prop:Beta_Tail_k}, we need another result which we prove for $k=2$ only:
\begin{lemma}\label{lem:Bd_ren_event}
   For $k =2$ there exists $q>0$ depending on $\kappa, \lambda$ such that for any $j \ge 1$, we have
   \begin{align}\label{eq:ren_ev}
       \bbP( \mathrm{Ren}_j \mid \mathbf{1}_{\mathrm{Ren}_1}, \mathbf{1}_{\mathrm{Ren}_2}, \ldots, \mathbf{1}_{\mathrm{Ren}_{j-1}}, \cF_{\tau_j}) = \bbP( \mathrm{Ren}_j \mid \cF_{\tau_j}) \ge q.
   \end{align}
\end{lemma}

\begin{proof}
 We first show that \eqref{eq:ren_ev} holds for $j=1$. We note that $\cF_{\tau_1}$ does not contain any 
 information about the PPP in the region $\mathbb{H}^+(W^{\text{move}}_{\tau_1}(2) + \kappa)$. For a single path, we observe
 \begin{align*}
     \tilde{q} &= \bbP \Big( \# \bigl ( B^+(\mathbf{0},1) \cap C_{\pi/4}(\mathbf{0}) \cap \mathcal{N}_{\tau_1 +1} \bigr ) = 
     \# \bigl ( B^+((0,-\kappa), \kappa + 1)  \cap \mathcal{N}_{\tau_1 +1} \bigr ) = 1 \Big) \\ &=  \lambda (\kappa +1) \exp({-\lambda (\kappa +1)}) < 1.
 \end{align*}
For two paths, however, $(\tilde{q})^2$ does not give an automatic lower bound since the paths may come close to each other, causing corresponding regions to overlap.

For $k=2$, we observe that the probability $\bbP(\mathrm{Ren}_1)$ is continuous in the locations of the projected points $(g_1^{\uparrow,1}, g_1^{\uparrow,2})$. Since $g_1^{\uparrow,1}(2)=g_1^{\uparrow,2}(2)$, the translation invariance nature of the model implies that $\bbP(\mathrm{Ren}_1)$ depends {\it continuously} only on $r := |g_1^{\uparrow, 1}(1) - g_1^{\uparrow, 2}(1)|$,
the distance
between the first coordinates of the projected points.
If $r>2(\kappa+1)$, then the regions $B^+(g_1^{\downarrow,1},\kappa+1)$ and $B^+(g_1^{\downarrow,2},\kappa+1)$ are disjoint. By the independence of the PPP on disjoint regions, for any $r>2(\kappa+1)$ we can take $q = (\tilde{q})^2 >0$. On the other hand, for any 
$r \in [0, 2(\kappa +1)]$ we observe that the probability of renewal event is strictly positive. Therefore, we set $q$ as
\begin{align*}
    q := \min \{ \bbP( \mathrm{Ren}_j \mid  |g_1^{\uparrow, 1}(1) - g_1^{\uparrow, 2}(1)|=r) : r \in [0, 2(\kappa +1)]\} \wedge (\tilde{q})^2 >0,
\end{align*}
to yield a uniform positive lower bound. This completes the proof for $j=1$.

 Next consider \eqref{eq:ren_ev} for $j=2$. Observe that the 
 r.v. $\mathbf{1}_{\mathrm{Ren}_1}$ depends only on the PPP
 in the half-plane $\mathbb{H}^-(W^{\text{move}}_{\tau_1}(2) + \kappa +1)$ and by definition we have $W^{\text{move}}_{\tau_2}(2) - W^{\text{move}}_{\tau_1}(2) > \kappa + 1$. Conditioned on $\cF_{\tau_2}$, since the subsequent steps as well as the 
 event $\mathrm{Ren}_2$ depend only on the PPP in the region
 $\mathbb{H}^+(W^{\text{move}}_{\tau_2}(2))$ we have
 \begin{align*}
     \bbP( \mathrm{Ren}_2 \mid \cF_{\tau_2}, \mathbf{1}_{\mathrm{Ren}_1}) = \bbP( \mathrm{Ren}_2 \mid \cF_{\tau_2}).
 \end{align*}
 Then the proof follows from similar argument as in the case
 of $j=1$. Finally, for general $j \ge 1$ the proof follows by
 method of induction.
\end{proof}

\noindent {\bf Proof of Proposition~\ref{prop:Beta_Tail_k}:}
We claim that for any $j \ge 1$ and $\ell \ge 0$
 \begin{equation}\label{eq:stoch_dom}
     \bbP( \gamma_{\ell +1} - \gamma_\ell > j \mid \mathcal{G}_\ell) \le \bbP( \mathbf{G} >j),
 \end{equation}
 where $\mathbf{G}$ is a geometric r.v. with success probability $q$. This means that conditionally on the $\sigma$-algebra $\mathcal{G}_\ell$, the random variable $\gamma_{\ell +1} - \gamma_\ell $ is
 stochastically dominated by $\mathbf{G}$. We prove the claim
 \eqref{eq:stoch_dom}
 for only $\ell =0$ case. The argument for general $\ell \ge 0$ is the same. Now observe that
 \begin{align*}
    \bbP( \gamma_1 - \gamma_0 > j \mid \mathcal{G}_0) &= \bbP( \mathrm{Ren}^c_1, \mathrm{Ren}^c_2, \ldots, \mathrm{Ren}^c_j \mid \mathcal{G}_0) \\
    &=\bbE\Big[\prod_{i = 1}^{j-1}\mathbf{1}_{\mathrm{Ren}^c_i} \bbE[\mathbf{1}_{\mathrm{Ren}^c_j} \mid \mathcal{S}_{j-1}, \cF_{\tau_j}] \Bigm| \mathcal{G}_0 \Big] \\
    &= \bbE\Big[\prod_{i = 1}^{j-1}\mathbf{1}_{\mathrm{Ren}^c_i} \bbE[\mathbf{1}_{\mathrm{Ren}^c_j} \mid \cF_{\tau_j}] \Bigm| \mathcal{G}_0 \Big] \\
    & \le (1-q) \bbE\Big[\prod_{i = 1}^{j-1}\mathbf{1}_{\mathrm{Ren}^c_i} \Bigm | \mathcal{G}_0 \Big] \le (1 - q)^j.
 \end{align*}
 We obtain the above bound by applying Lemma~\ref{lem:Bd_ren_event} and
 a use of recursion. We observe that conditioned on  the $\sigma$-algebra $\cS_j$, the difference $\tau_{j+1} - \tau_j$ has exponentially decaying tail. This follows from Proposition~\ref{prop:Tau_tail_k} and the observation that the events $\mathrm{Ren}_1, \ldots, \mathrm{Ren}_{j-1}$ depend only on the
 Poisson points in $\mathbb{H}^-(W^{\text{move}}_{\tau_j}(2) + \kappa +1)$ while $\tau_{j+1} - \tau_j$ depends on Poisson points in $\mathbb{H}^+(W^{\text{move}}_{\tau_j}(2))$ giving:
 \begin{align*}
     \bbP(\tau_{j+1} - \tau_j >n \mid \mathbf{1}_{\mathrm{Ren}_1}, \ldots, \mathbf{1}_{\mathrm{Ren}_{j-1}}, \cF_{\tau_j} ) =  \bbP(\tau_{j+1} - \tau_j >n \mid \cF_{\tau_j} ) \le c_0 \exp(-c_1 n).  
  \end{align*}
We highlight that the decay constants $c_0, c_1 > 0$ depend only on $\lambda, \kappa$. 
 
 Now, we can construct a random variable $\Delta$ such that $\bbP( \Delta >n) \le  c_0 \exp(-c_1 n) $ and stochastically 
 dominates the conditional difference $\tau_{j+1} - \tau_j >n \mid \cS_j$. Moreover, we can stochastically bound, for any
 $\ell, m$ the sum
 \begin{align*}
     \sum_{j=0}^{m-1} \tau_{\gamma_\ell + j + 1} - \tau_{\gamma_\ell + j }
 \end{align*}
 conditionally on $\cG_\ell$ by $\sum_{i=1}^m \Delta_i$  where  $\Delta_i$'s are i.i.d. copies of $\Delta$. Now choose  $\eta>0$ small enough so that $\bbE[\exp(\eta \Delta )] < \infty$. Then for any constant $c>0$, we have
 \begin{align*}
     &\bbP( \beta_{\ell +1} - \beta_\ell >n \mid \cG_\ell) \\
     & \le \bbP\Big(\sum_{j=0}^{\lfloor cn \rfloor} \tau_{\gamma_\ell + j + 1} - \tau_{\gamma_\ell + j } > n \Bigm| \cG_\ell \Big) + \bbP( \mathbf{G} > \lfloor cn \rfloor \mid \cG_\ell) \\
     & \le \bbP \Big( \sum_{j=1}^{\lfloor cn \rfloor}  \Delta_j > n \Big) + ( 1 - q )^{\lfloor cn \rfloor} \\
     & \le e^{-\eta n} \bbE[\exp(\eta \Delta)]^{\lfloor cn \rfloor}  + ( 1 - q )^{\lfloor cn \rfloor}. 
 \end{align*}
 Hence, the proof of \eqref{eq:beta_l} is complete by taking
 $c= c(\eta)$ sufficiently small.
\qed

Let $\mathbf{x}^1, \ldots, \mathbf{x}^k$ be starting points of $k$ trajectories and fix
$\ell \ge 0$. We want to understand the behavior of random rectangles containing the
regions explored by the $k$($=2$) trajectories of the joint exploration process $\{g_n(\mathbf{x}^1), \ldots, g_n(\mathbf{x}^k) : n \ge 0 \}$ between the $\ell$- th and $(\ell + 1)$-th renewal steps. To this end we define,
\begin{align}\label{def:ren_block_size}
    \mathbf{W}_{\ell +1} =  \mathbf{W}_\ell( \mathbf{x}^1, \ldots, \mathbf{x}^k) := \sum_{m= \beta_\ell}^{\beta_{\ell +1} - 1} \|W^{\text{move}}_m - h(W^{\text{move}}_m) \|.
\end{align}
By definition the r.v. $ \mathbf{W}_{\ell +1} $ is such that for any $1 \le i \le k$, the 
random set
\begin{align*}
    \bigcup_{m= \beta_\ell}^{\beta_{\ell + 1} - 1} B^+( g_m(\mathbf{x}^i), \| g_m(\mathbf{x}^i)- g_{m+1}(\mathbf{x}^i) \|),
\end{align*}
is included in the region $g_{\beta_\ell} (\mathbf{x}^i) \circ ([-\mathbf{W}_{\ell +1}, \mathbf{W}_{\ell +1}] \times [0, \mathbf{W}_{\ell +1}] )$. The r.v. $ \mathbf{W}_{\ell +1} $ dominates \textit{width} of the renewal
blocks. In summary, 
the event 
$$
\{ \cap_{i=1}^k  (g_{\beta_\ell} (\mathbf{x}^i) \circ ([-\mathbf{W}_{\ell +1}, \mathbf{W}_{\ell +1}] \times [0, \mathbf{W}_{\ell +1}])) = \emptyset \}
$$
implies that in between the $\ell$-th and $(\ell + 1)$-th joint renewals, the two DSF paths evolve independently. 

The same argument as in Lemma~\ref{lem:exp_tail_width} gives us the next result which shows that the conditional distribution of $\mathbf{W}_{\ell + 1} \mid \cG_\ell$ is stochastically dominated by a random variable $\mathrm{W}$ which has
sub-exponential tail decay.

\begin{proposition}\label{prop:exp_decay_size}
    There exists constants $c_0, c_1 >0$ depending on $\lambda, \kappa$ and $k$ such that for any 
    $\ell \ge 0$ and $n \in \bbN$, we have
    \begin{equation}
        \bbP( \mathbf{W}_{\ell + 1} > n \mid \cG_\ell ) \le c_0 \exp(-c_1 n^{1/2}).
    \end{equation}
\end{proposition}

If DSF paths are sufficiently far apart at the $\beta_\ell$-th step, then the occurrence of the renewal event ensures that we can continue with the restarted DSF paths starting from the projected points $g^{\uparrow,1}_{\beta_\ell}, \cdots, g^{\uparrow,k}_{\beta_\ell}$ conditional that all of their next steps are top steps within respective $\ell_\infty$ balls of radius $1$ each. On the other hand, at the $\beta_\ell$-th step if the DSF paths are closer, then future evolution of DSF paths may not match with that of the restarted paths. Still in that case, as observed in \cite{MR4203342} in case of $\ell_2$ DSF paths, the coalescing time of the restarted paths would dominate the coalescing time of the actual DSF paths. Therefore, we would work with the process of restarted paths at renewal steps and towards that we define for $\ell \ge 0$, 
\begin{align}
    \label{def:Z_ell}
    Z_\ell := g^{\uparrow, 2}_{\beta_\ell}(1) - g^{\uparrow, 1}_{\beta_\ell}(1).
\end{align}
The same argument as Corrolary 4.8 in \cite{MR4203342} proves the following lemma which says that, far from the origin, 
the process $\{Z_\ell : \ell \geq 0\}$ behaves like a mean zero 
random walk satisfying certain moment bounds. 

\begin{lemma}
\label{lem:ZprocessRwalkProperties}
Fix $\mathbf{x}, \mathbf{y} \in \mathbb{R}^2$ with $\mathbf{x}(2) =  \mathbf{y}(2)$, $\mathbf{x}(1) <  \mathbf{y}(1)$. Consider the joint exploration process of $\ell_\infty$ DSF paths starting from $\mathbf{x}, \mathbf{y}$ till the $\ell$-th (joint) renewal step. 
Given the $\sigma$-algebra ${\cG}_{\ell}$, 
there exist positive constants $M_0, C_0, C_1, C_2$ and $C_3$ 
and an event $\mathbf{F}_\ell$ such that:
\begin{itemize}
\item[(i)] On the event $\{Z_\ell > M_0 \}$ we have 
$$
\mathbb{P}( \mathbf{F}^c_\ell \mid {\cG}_{\ell} ) \leq C_3/(Z_\ell)^3\text{ and }
\mathbb{E} \big[ (Z_{\ell +1} - Z_{\ell}) \mathbf{1}_{\mathbf{F}_\ell} \mid {\cG}_{\ell} \big] = 0.
$$
\item[(ii)] On the event $\{ Z_\ell \leq  M_0 \}$ we have
 $$
 \mathbb{E} \big[ | Z_{\ell +1} - Z_{\ell} |  \mid {\cG}_\ell \big] \leq C_0 .
 $$

\item[(iii)] For any $\ell\geq 0$ and $m>0$, there exists $c_m \in (0,1)$ such that, 
on the event $\{Z_\ell \leq m \}$,
\begin{equation*}
\mathbb{P} \big( Z_{\ell+1} = 0 \mid {\cG}_\ell \big) \geq c_m ~.
\end{equation*}
\item[(iv)] On the event $ \{Z_\ell > M_0 \}$, we have
\begin{equation*}
\mathbb{E} \bigl[ ( Z_{\ell+1} - Z_{\ell} )^2 \mid {\cG}_\ell \bigr] \geq C_1 \; \text{ and } \;
\mathbb{E} \bigl[ |Z_{\ell+1} - Z_{\ell} |^3 \mid {\cG}_\ell  \bigr] \leq C_2 ~.
\end{equation*}
\end{itemize}
\end{lemma}

\begin{proof}
Take $M_0$ large and recall the definition of the random variable $\mathbf{W}_\ell$ from \eqref{def:ren_block_size}. We consider the event $ \mathbf{F}_\ell = \{ \mathbf{W}_{\ell +1} \le Z_\ell /3 \} $.  Proposition~\ref{prop:exp_decay_size} gives us that $\mathbb{P}( \mathbf{F}^c_\ell \mid {\cG}_{\ell} ) \leq C_3/(Z_\ell)^3$, for some constant $C_3 >0$.  

Now on the event $\{ Z_\ell \ge M_0, \mathbf{W}_{\ell +1} \le Z_\ell /3 \}$ consider the trajectories of the two paths starting from $\mathbf{x}^1$ and $\mathbf{x}^2$ between the $\ell$-th and $(\ell+1)$-th renewal steps. These regenerated trajectories, till the next (joint) renewal step, can be constructed using PPP on the region $\mathbb{H}^+( W^{\text{move}}_{\beta_\ell}(2) + \kappa)  $ within disjoint rectangles 
$$
g_{\beta_\ell}^{\uparrow, 1} \circ ([-\mathbf{W}_{\ell + 1}, \mathbf{W}_{\ell + 1}] \times [0, \mathbf{W}_{\ell + 1}])\text{ and }
g_{\beta_\ell}^{\uparrow, 2} \circ ([-\mathbf{W}_{\ell + 1}, \mathbf{W}_{\ell + 1}] \times [0, \mathbf{W}_{\ell + 1}]) .
$$ 
We construct a new PPP $\cN^*$ as follows:
\begin{itemize}
    \item[a)] Conditioned on $\cG_\ell$, we interchange the realizations of the PPP in $\mathbf{R}_1 = g_{\beta_\ell}^{\uparrow, 1} \circ ([- Z_\ell/3, Z_\ell /3] \times [0, Z_\ell /3])$ and  
    $\mathbf{R}_2 = g_{\beta_\ell}^{\uparrow, 2} \circ ([- Z_\ell/3, Z_\ell /3] \times [0, Z_\ell /3])$. 
    \item[b)] We keep the realization of $\cN$ outside $\mathbf{R}_1\cup\mathbf{R}_2$ unchanged.
\end{itemize}
We observe that for the newly constructed PPP $\cN^*$, each of the rectangles $\mathbf{R}_1$ and $\mathbf{R}_2$ contains the unique Poisson point that is uniformly distributed in $B^+(g_{\beta_\ell}^{\uparrow,i},1)$ for $i=1,2$ respectively. Consequently, $\cN^*$ has the same distribution as a PPP conditioned to contain a unique point in each semi-rectangle $B^+(g_{\beta_\ell}^{\uparrow,i},1)$ and no points in $\mathbb{H}^+( W^{\text{move}}_{\beta_\ell}(2) + \kappa) \cap B^+(g_{\beta_\ell}^{\downarrow, i}, \kappa +1)) \setminus B^+(g_{\beta_\ell}^{\uparrow, i}, 1)$ for $i=1,2$.

To this end, consider the event $\mathbf{F}_\ell$ and construct the trajectories between the $\ell$-th and $(\ell+1)$-th renewal steps using $\cN^*$. The number of steps until the next renewal time and the size of the corresponding renewal block remain unchanged. However, the increments of the two paths over this renewal interval are interchanged. Therefore, the increment $Z_{\ell+1}-Z_\ell$ is replaced by 
$-(Z_{\ell +1 }- Z_\ell)$, which completes the proof.

   Item $(ii)$ instantly follows from the following observation
   \begin{align*}
       \bbE[ Z_{\ell +1} - Z_\ell \mid \cG_\ell] \le \bbE[| Z_{\ell +1} - Z_\ell | \mid \cG_\ell] \le \bbE[ 2 \mathbf{W}_{\ell + 1} \mid \cG_\ell] < \infty,
   \end{align*}
since conditioned on $\cG_\ell$,  $\mathbf{W}_{\ell + 1}$ admits sub-exponential tail
decay, see Proposition~\ref{prop:exp_decay_size}.

For item $(iii)$, we note that the $\sigma$-algebra $\cG_\ell$ does not contain any 
information about the region $\mathbb{H}^+( W^{\text{move}}_{\beta_\ell}(2) + \kappa) \setminus (B^+(g_{\beta_\ell}^{\downarrow, 1}, \kappa +1) \cup B^+(g_{\beta_\ell}^{\downarrow, 2}, \kappa +1)) $. Thus it is not difficult to see that 
the conditional probability $\bbP( Z_\ell \mid \cG_\ell)$ is strictly positive.

   Finally, it remains to prove item $(iv)$. We note that
   \begin{equation*}
       \bbE[ (Z_{\ell +1} - Z_\ell)^2 \mid \cG_\ell] \ge  \bbE[ (Z_{\ell +1} - Z_\ell)^2 \mathbf{1}_{\{|Z_{\ell +1} - Z_\ell| > 1 \}} \mid \cG_\ell] \ge \bbP ( |Z_{\ell +1} - Z_\ell|^2 > 1 | \cG_\ell).
   \end{equation*}
 Further, one can easily check that on the event $\{Z_\ell > M_0\} $, the probability 
 $\bbP((Z_{\ell +1} - Z_\ell)^2 \mid \cG_\ell)$ is strictly positive. For the third moment,
 \begin{equation*}
    \bbE[ (Z_{\ell +1} - Z_\ell)^3 \mid \cG_\ell] \le \bbE[| Z_{\ell +1} - Z_\ell |^3 \mid \cG_\ell] \le \bbE[ (2 \mathbf{W}_{\ell + 1})^3 \mid \cG_\ell] < \infty. 
 \end{equation*}
\end{proof}

Lemma \ref{lem:ZprocessRwalkProperties} enables us to prove Theorem \ref{thm:DSF_CoalTime_Decay}.

\vspace{0.5cm}

\noindent{\bf Proof of Theorem \ref{thm:DSF_CoalTime_Decay}:}
Lemma \ref{lem:ZprocessRwalkProperties} allows us to apply Corollary 5.6 in \cite{MR4203342} for two $\ell_\infty$ DSF paths starting from $\mathbf{x}, \mathbf{y} \in \mathbb{R}^2$ with
$\mathbf{x}(2) = \mathbf{y}(2)$ and gives the required decay estimate on the coalescing time tail. This completes the proof. 
\qed

\vspace{0.5cm}

Now we are left with proving Theorem~\ref{prop:DSF_Tree}.

\vspace{0.5cm}

\noindent {\bf Proof of Theorem~\ref{prop:DSF_Tree}:}
Recall that from Theorem~\ref{thm:DSF_CoalTime_Decay}, we know
that for all $\mathbf{x}, \mathbf{y} \in \bbR^2$ with $\mathbf{x}(2) = \mathbf{y}(2)$ the two trajectories originated from $\mathbf{x}$ and $\mathbf{y}$ coalesces almost surely. Let $\mathbf{u}, \mathbf{v} \in \cN$ and fix $t >0$
such that $t > \max \{\mathbf{u}(2), \mathbf{v}(2) \}$. Since,
$t > \mathbf{u}(2)$ and $h^n(\mathbf{u})(2) \to \infty$ as 
$n \to \infty$, there exists a unique integer $n_1 \ge 1$ such that
\begin{align*}
    h^{n_1}(\mathbf{u})(2) \le t \text{ and }  h^{n_1 +1}(\mathbf{u})(2) > t.
\end{align*}
By the construction of DSF, there exists a Poisson Point
$\mathbf{x} \in \cN$ with $\mathbf{x}(2) = t$ whose first
ancestor is $h^{n_1 +1}(\mathbf{u})$, i.e., $h(\mathbf{x}) = h^{n_1 +1}(\mathbf{u})$ almost surely. Similarly, there
exists an integer $n_2 \ge 1$ and a Poisson Point $\mathbf{y} \in \cN$ with $\mathbf{y}(2) = t$ satisfying $h(\mathbf{y}) = h^{n_2}(\mathbf{v})$. Since $\mathbf{x}$ and $\mathbf{y}$ lie on the same horizontal level, Theorem~\ref{thm:DSF_CoalTime_Decay} implies that the trajectories originating from $\mathbf{x}$ and $\mathbf{y}$ coalesce almost surely. As these trajectories coincide with the tails of the trajectories originating from $\mathbf{u}$ and $\mathbf{v}$, respectively, it follows that the trajectories from $\mathbf{u}$ and $\mathbf{v}$ also coalesce almost surely. This completes the proof. 
\qed

\section*{Acknowledgements}
Dipranjan Pal expresses his gratitude to Indian Statistical Institute for providing a fellowship to support his PhD studies. Kumarjit Saha gratefully acknowledges
financial support via the CEFIPRA grant 6901/1.

\bibliographystyle{amsplain}
\bibliography{example}

\end{document}